\documentclass[11pt]{article}

\usepackage[margin=1in]{geometry}
\usepackage{amsmath,amsthm,amssymb}
\usepackage{commath,mathtools,bm}
\usepackage{algorithm,algorithmic}
\usepackage{graphicx}
\usepackage{subcaption}
\usepackage{hyperref}
\usepackage{color}
\usepackage[sort,nocompress]{cite}

\graphicspath{ {./images/} }

\title{Approximating High-Dimensional Minimal Surfaces with Physics-Informed Neural Networks\thanks{SZ is currently a 12th grade student at GSMST.}}
\author{Steven Zhou\footnote{Gwinnett School of Mathematics, Science, and Technology (GSMST), Lawrenceville, GA 30044. Email: \url{steven.zhou1313@outlook.com}} \and Xiaojing Ye\footnote{Department of Mathematics and Statistics, Georgia State University (GSU), Atlanta, GA 30303. Email: \url{xye@gsu.edu}}}

\begin{document}

\maketitle

\begin{abstract}

In this paper, we compute numerical approximations of the minimal surfaces, an essential type of Partial Differential Equation (PDE), in higher dimensions. Classical methods cannot handle it in this case because of the “Curse of Dimensionality”, where the computational cost of these methods increases exponentially fast in response to higher problem dimensions, far beyond the computing capacity of any modern supercomputers. Only in the past few years have machine learning researchers been able to mitigate this problem. The solution method chosen here is a model known as a Physics-Informed Neural Network (PINN) which trains a deep neural network (DNN) to solve the minimal surface PDE. It can be scaled up into higher dimensions and trained relatively quickly even on a laptop with no GPU. Due to the inability to view the high-dimension output, our data is presented as snippets of a higher-dimension shape with enough fixed axes so that it is viewable with 3-D graphs. Not only will the functionality of this method be tested, but we will also explore potential limitations in the method's performance.
\end{abstract}
\large

\section{Introduction}Soap bubbles are the most common examples of minimal surfaces in our everyday lives,
but have you ever wondered what they would look like if the frame creating the bubble was
in 3 or more dimensions? Minimal surfaces can be visualized easily in our three-dimensional space. Mathematically, we can set 2-dimensional boundary conditions and use equations to evaluate each point in said 2-D domain and add a third dimension in a way respecting the rules governing what counts as a minimal surface. But directly visualizing minimal surfaces above 3 dimensions is impossible. Taking 3-dimensional figures as a boundary for extruding our minimal surface into the 4th dimension is simply unfathomable for the human mind, let alone even more dimensions. Fortunately, with the power of modern computing, we can solve these minimal surfaces regardless of our human limitations for a small first step into visualizing the higher dimensions.

Partial Differential Equations (PDE) have been used to represent complex systems involving many physics-related curiosities, including minimal surfaces \cite{berg2018unified}. They are one of the fundamental building blocks of what makes each minimal surface truly ``minimal'', mathematically defining how the total surface area is always the smallest value possible. Using them we can consistently build minimal surfaces for every situation, but the possibilities of simply solving PDEs are limited by certain factors. 

The ``Curse of Dimensionality'', a plague on all computing solutions of high-dimension PDEs consistently limits our computing models. As the dimensions of the data increase, so do the required amount of data points. This exponential growth of these requirements means that higher-dimension equations will require astronomical amounts of computing power in order for the calculations to remain accurate. Otherwise, all results will become sparse -- insufficient to provide a consistent solution across the entire domain. There are simply too many data points needed for us to continue on without new options.  

Machine learning has become a prevalent resource for the computation of complex mathematical figures. Neural networks (NN), a specific type of function approximation model, are one of the great tools used by modern mathematicians to approximate an equation with no analytical solution. Multi-dimensional inputs and loss function minimization techniques allow us to handle equations with many variables without fear of wasting valuable computational power. 

Earlier versions of using Deep Neural Networks to solve PDEs have and could only place them in the roles of assisting the conventional methods of finite elements and finite differences \cite{nabian2018deep}. A lack of data to train an effective neural network model is still one of the limiting factors to traditional deep-learning PDE solutions when utilized in the real world. But in the past few years, even the issue of a lack of comprehensive data has been mitigated by the development of Physics Informed Neural Networks (PINN), which track the physical attributes of the data \cite{raissi2019physics-informed, raissi2017physics, hu2022general}. PINN in general requires a multitude of conditions: initial values for the solution, governing equations (so that the information reflects a continuous model), and values that are of importance on the boundaries \cite{pang2020npinns}. Famous PDEs have been approximated such as the Navier-Stokes equation \cite{deryck2023stokes}. Also, PINN is a significantly simpler species of physics-informed learning program than similar algorithms such as physics-informed Gaussian processes (PIGI). These specific Neural Networks and their various formats, including Sparse Physics-based partially-Interpretable Neural Networks (SPINN) \cite{ramabathiran2021spinn} and Fractional Physics Informed Nerual Networks (fPINN) \cite{pang2019fpinns} retain gradients over many iterations of training and observe the patterns between the results of data with similar value. The creators of SPINN, seeking a type of neural network that would be easy to observe the inner workings of, created a model with fewer inner connections between neural network layers. Meanwhile, fPINN allows a seemingly tiny amount of data at certain points to be used to create a representation of the full equation. Both methods can be used in a wide array of scientific pursuits but are overly qualified for the results we wish to achieve. A few other examples of recently-developed PINN methods include Conservative Physics-Informed Neural Networks (cPINN) \cite{jagtap2020conservative}, Extended Physics-Informed Neural Networks (XPINN) \cite{jagtap2020xpinn, jagtap2022supersonic, shukla2021parallel, penwarden2023sweeping, mishra2021pde}, and Augmented Physics-informed Neural Networks (APINN) \cite{hu2022augmented}. Each of these new formats of PINN is associated with different classes of PDEs as a specialized approach to solving them. 

PINN has also been widely applied to fields such as analyzing micro-structures \cite{shukla2022micro}, fluid mechanics
\cite{jagtap2022waves}, and aerodynamic flow
\cite{mao2020highspeed}.While PINN has been used to track equations in \cite{raissi2019physics-informed}, the authors did not use it to solve the minimal surface equation. At the time of writing this paper, \cite{kabasi2023membranes} had been the only published paper using PINN for calculating 2D domain minimal surfaces, but here we will apply PINN to compute and visualize minimal surfaces in higher dimensions. We wish to transcend the boundaries of what PINN has done so far and contribute to the movement in which machine learning has emerged as one of the main tools for solving high-dimension PDEs, with minimal surfaces being just another example.

With the ability of PINN to create a surface with logical definitions of the data's behavior at every point, we can use a multi-layer neural network based on the minimal surface equation \cite{magill2018neural} to model the desired equation into what is often visually represented as a soap bubble, with two dimensions for the domain and given values on its boundaries. A unique view of soap bubbles in higher dimensions can be formed by beginning with a 1 by 1 square and adding in unit dimensions as the domain increases to 3D and 4D. The boundary conditions on the sides of the square in 2D and cube in 3D act like the wire frames of a bubble mold, confining the resulting minimal surface(the bubble) to a certain region. 

Our last step is to visualize the newly-computed minimal surfaces. Since we cannot see into 4-dimension space, we will take "slices" of the bubble from the faces of each minimal surface where we can assume one or two variables to be constant. This method will be employed in our 3D and 4D examples. And due to the sensitive nature of training neural network approximations, our work in this paper is heavily based on constantly adjusting the weights for the two terms of the loss function, namely the 
square norm of residual for the equation and the discrepancy in boundary fitting. Also, the limitations of using PINN will be explored in detail, especially in situations of higher dimensions skewing surface generation.  



Machine learning methods such as the Deep Ritz method \cite{e2018deep} and  weak adversarial networks  \cite{zang2020weak, bao2020numerical} can also be used by lowering the computational power required for PDEs and create a roundabout method for accurately approximating solutions. Both are feasible alternatives for PINN in this scenario, but as PINN has more wide-spread applications, it has become our method of choice.

The rest of the paper is as follows. Section 2 will define neural networks and explain the specifics of how they approximate functions. We will also introduce the specific neural network architecture we use and how it is layered. In section 3, we construct the boundary frame in 2D, 3D, and 4D. We also discuss the loss function that gives the neural network the namesake of physics-informed and the method for visualizing minimal surfaces created with 3D and 4D domains in a 2D domain. Furthermore we discuss the steps taken to implement the PINN. Section 4 discusses the results and the accuracy of said examples in 3 different dimension inputs. It devotes effort to showing that a specific mixture of weights for each boundary equation can produce mostly optimal results with negligible error.  Examples are provided to prove our model retains the ability to calculate minimal surfaces even in the presence of clashing boundary inputs. Also included are scenarios where the model breaks down. Finally, Section 5 concludes the paper with a brief summary. 

\section{Neural Networks}

To understand the usage of PINN, it is important to quickly summarize the notion of Neural Networks (NN). In every neural network, information is first input in the form of a vector, which through multiple layers of transmutation between the different neurons, eventually reach an output layer approximating the results if said input had been run through a complex mathematical function. In a way, NN acts like a Taylor polynomial, approximating a function far more complicated than one possibly writable or even fathomable with the tools of mathematics.

We can model the deep neural network as a machine learning algorithm consisting of $L$ layers in the form 
\begin{equation}\label{eq:dnn}
y = W_{L}h_{L-1} + b_{L}, 
\end{equation}
where the input layer is $h_0 = x$ and the $l$th hidden layer is
\begin{equation*}
    h_l= \sigma (W_{l}\, h_{l-1} + b_{l})
\end{equation*}
for $l=0,1,\cdots,L-2$.
In equation \eqref{eq:dnn}, $L$ represents the total layer number of the DNN. The result $h_l$ is called the output of $l$th hidden layer, which will be used as the input to compute $h_{l+1}$ for the next, $(l+1)$th layer. In this format, $W_{l}$ plays the role of a vector dimension converter, converting the dimension of the vector $h_{l-1}$ to that of $h_l$. $b_{l}$, known as the bias, is a real valued vector the same size of $h_l$ which shifts the values in each hidden layer.
The nonlinear activation function $\sigma$ is used to filter the result $W_l h_{l-1} + b_l$ in each layer to obtain the output $h_l$. This function must be nonlinear and is often applied component-wise. Some of the most famous activation functions include Sigmoid, hyperbolic tangent (Tanh), Rectified Linear Unit (ReLU), and Exponential Linear Unit (ELU).  Hyperbolic tangent comes in the form of $\frac{\sinh{x}}{\cosh{x}}$. Even more simplified, the activation function can be taken to equal $\frac{e^{2x}-1}{e^{2x}+1}$, resulting in most positive terms being replaced with a value around 1 while negative terms are set to a value close to -1. Rectified Linear Unit, on the other hand, raises all negative values within the set to 0 and does not affect the nonnegative values in said vector set. 

The $(l+1)$th layer in the DNN takes input from a certain number of ``neurons'', represented by $h_l$, and transmutes it across a chosen amount of new neurons so that over the course of several layers, all possibilities are accounted for and the network will produce a balanced result for any desired input. At the beginning of this process, sets of multi-variable inputs (2-4 variables in our cases) will be slowly widened to a total of 128 neurons per set before being condensed back to the output value $u(x;\theta)$ consisting of only one variable for each original set of input variables. Here $\theta$ is the collection of parameters including $W_{l}$ and $b_{l}$ for all $l=1,\dots,L$, a representation that the solution we find is only an approximation of $u$, not the exact value.

\begin{figure}
\centering
\includegraphics[width = 0.6\textwidth]
{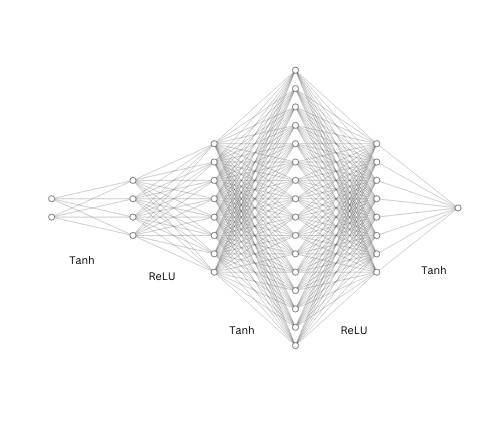}
\caption{A simplified diagram of the Neural Network structure used. The input layer is variable between 2, 3, and 4 nodes; the hidden layers contain in order 32, 64, 128, and 64 nodes. Tanh and ReLU activation functions are used between layers. This architecture is selected as an efficient model for PINN.}
\label{neuralnetwork}
\end{figure}

After one iteration through the multiple layers of the neural network, the model must be trained by an optimizer before the next epoch of testing. In every neural network, an optimizer and loss function work together to measure the shortcomings of each layer in the neural network and accordingly adjust the values of $W_{l}$ and $b_{l}$ for all $l$. The term ``backpropagation'' is used to refer to this process of finding the gradients of each layer from the output until the input \cite{nabian2018deep}. In our current paper we will instead use the Adam optimizer that yields fast results \cite{kingma2014adam:}. 

\section{Methods and Model}
\subsection{Minimal Surfaces}

To calculate our minimal surface, we can first define it in a minimization framework, 
\begin{equation}    
\min_{u} \cbr[2]{ E(u) := \int_{\Omega} \sqrt{1 + |\nabla u(x)|^2} \, dx\  \quad u(x)  = g(x), \;\; x \in \partial \Omega}
\label{eq:minimal}
\end{equation}
where $u(x)$ represents the minimal surface and $\nabla u$ is the gradient of $u$. $\partial\Omega$ refers to the boundary lines of the domain where one of the dimensions has a value of 0 or 1 in our case. $g(x)$ is a given function defined on the boundary representing the frames of a soap bubble. The solution $u(x)$ can be computed through the following equation, 
\begin{align}
    -\nabla \cdot \left(\frac{\nabla u(x)}{\sqrt{1 + |\nabla u(x)|^2}}\right) & = f(x), \quad x \in \Omega \label{eq:laplace}\\
    u(x) & = g(x), \quad  x \in \partial \Omega \label{eq:dirichelt_bdry}
\end{align}
Here $\Omega \subset \mathbb{R}^d$ is a bounded domain in however many input dimensions given, and $-\nabla \cdot$ is the negative divergence. For the sake of simplicity, we will begin with a 2D representation in the form of a square. In general, the ``frame'' will be a square curving into the 3rd dimension for our 2D frame and a cube with each edge extruding into the 4th dimension for the 3D model. Both $f(x)$ and $g(x)$ are known given functions. In equation \eqref{eq:minimal}, $f(x)$ is set to 0 for all minimal surfaces in an Euclidean space.

\subsection{Physics-Informed Neural Networks}

The equations \eqref{eq:laplace} and \eqref{eq:dirichelt_bdry} determine our minimal surface $u(x)$, but to approximate it with our PINN we must first reduce it to a usable form. 

We can get an equivalent set of equations
\begin{align}
        &\min_{u}\quad  \int_{\Omega} \Big| -\nabla \cdot \left(\frac{\nabla u(x)}{\sqrt{1 + |\nabla u(x)|^2}}\right) - f(x)\Big|^2 \dif x 
        \label{eq:firstpart}\\
        &\mbox{subject to}\quad u(x) = g(x), \quad \forall\, x \in \partial \Omega
        \label{eq:secondpart}
\end{align}
to minimize the interior error. Continuing on, we can also add equation \eqref{eq:secondpart} into the minimization as an additional penalty term to get

\begin{equation}
    \min_{u} \quad \int_{\Omega} \Big| -\nabla \cdot \left(\frac{\nabla u(x)}{\sqrt{1 + |\nabla u(x)|^2}}\right) - f(x)\Big|^2 \dif x + \alpha \int_{\partial \Omega}|u(x) - g(x)|^2 \dif S(x)
\label{eq:beforepinn}
\end{equation}
where $\alpha>0$ is a calibrated constant weight that we, the experimenters, will manipulate to balance the boundary loss and any internal loss.

Our PINN approximation for $u(x)$ can be parameterized as the neural network $u(x;\theta)$ in a space $\Omega$ to form the loss function for the minimal surface equation
\begin{equation}
    \label{eq:loss}
    \ell(\theta) = \tilde{L}(u(\cdot;\theta)) = \int_{\Omega} \Big| -\nabla \cdot \left(\frac{\nabla u(x;\theta)}{\sqrt{1 + |\nabla u(x;\theta)|^2}}\right) - f(x)\Big|^2 \dif x + \alpha \int_{\partial \Omega} | u(x;\theta) - g(x)|^2 \dif S(x).
\end{equation}

To train the neural network $u(\cdot;\theta)$, we need to find the minimizer $\ell(\theta)$ as the parameter of $u(\cdot;\theta)$. However, the integrals in equation \eqref{eq:loss} cannot be directly calculated by hand. Instead, we use a method known as Monte Carlo integration to evaluate these integrals. To that end, we can randomly sample a number of points from inside the square $\Omega$
and a number of points from the boundary $\partial\Omega$, and estimate the integrals using the averages of the integrands evaluated at these sampled points. 
More precisely, we have the following approximations:
\[
    \int_{\Omega} \Big| -\nabla \cdot \left(\frac{\nabla u(x;\theta)}{\sqrt{1 + |\nabla u(x;\theta)|^2}}\right) - f(x)\Big|^2 \dif x \approx \frac{1}{N} \sum_{i=1}^N \Big| -\nabla \cdot \left(\frac{\nabla u(x_i;\theta)}{\sqrt{1 + |\nabla u(x_i;\theta)|^2}}\right) - f(x_i)\Big|^2 
\]
for the internal loss where $x_i$'s are uniformly sampled in $\Omega$ and 
\[
    \alpha\int_{\partial\Omega} |u(x;\theta) - g(x)|^2 \dif S(x) \approx \frac{\alpha}{M} \sum_{i=1}^M | u(x_i; \theta) - g(x_i)|^2 
\]
for the boundary loss where $x_i$'s are uniformly sampled on $\partial \Omega$. These approximations can be viewed as an analog to the Mean Squared Error (MSE) losses used in typical deep learning problems.

To evaluate $u(x_i;\theta)$ we feed $x_i$ through the neural network $u(\cdot;\theta)$. To use the Adam optimizer to find the minimizer of $\ell(\theta)$, we can backpropagate the neural network and compute the gradient of $\ell$ in each iteration. This calculation can be done using automatic differentiation in modern deep learning frameworks such as PyTorch and TensorFlow. Through thousands of iterations of this process, the parameter $\theta$ gradually improves and the neural network $u(\cdot;\theta)$ will eventually get close to the solution $u(x)$.


\subsection{Challenges in implementation}

\begin{algorithm}
\begin{enumerate}
    \item Construct the Neural Network $u(x;\theta)$.  
    \item Input the boundary condition $g(x)$.
    \item Construct a frame in 2D, 3D, or 4D consisting of boundaries of $n$ points where one of the variables $x_i$ is 0 or 1 while the other(s) is a/are random real number(s) ranging from 0 to 1.
    \item Generate $m$ random points taken from the interior of the domain.
    \item Create a 1-dimension vector $f(x)$ with length $m$ where every entry is set to the same constant value.
    \item Create two copies of the boundary frame to run through $g(x)$ and the neural network separately
    \item Define loss function $\ell$ as in \eqref{eq:loss}. The gradient and divergence are implemented using automatic differentiation.
    \item Choose an optimizer, Adam in this scenario.
    \item \textbf{Gradient Descent: } Train the model using the Adam optimizer. 
    In each iteration of Adam, we use the loss equation's backward propagation based on automatic differentiation to compute the gradient of $\ell$ and update the optimizer accordingly. Repeat this for a given number of iterations. 
\end{enumerate}
\caption{PINN Computing Process}
\label{alg:pinn}
\end{algorithm}

To implement our Physics-Informed Neural Network, we can follow the steps in Algorithm \ref{alg:pinn}. This algorithm is a step-by-step pseudocode-like basis of how the neural network is created and trained.
Some of the steps must be modified between examples of different domain dimensions to account for changes such as the boundary frame shape.

All of these procedures can be performed on a GPU-less personal laptop quite easily within a few minutes. The general format of the gradient descent style is to take two sets of coordinates for the interior and exterior and create separate versions of each. The versions will represent the functions we wish to approximate and the current PINN approximation. By utilizing our loss function we can continually generate the loss gradients and train our neural network until it converges. 

A few issues when programming the neural network in Python arose regarding equation \eqref{eq:firstpart}. Calculating the norm in the dividend and applying it to each element of the gradient leads to consistently incorrect PDE loss values that did not converge. The problem eventually was resolved by re-configuring the norm calculations to only use values from a single data point for its corresponding gradient.

In order to prove the ability of the PINN algorithm, we have prepared a few examples.

\section{Results}

\begin{figure}
    \centering
    \begin{minipage}[b]{0.3\textwidth}
    \includegraphics[scale=0.4, width = \textwidth]{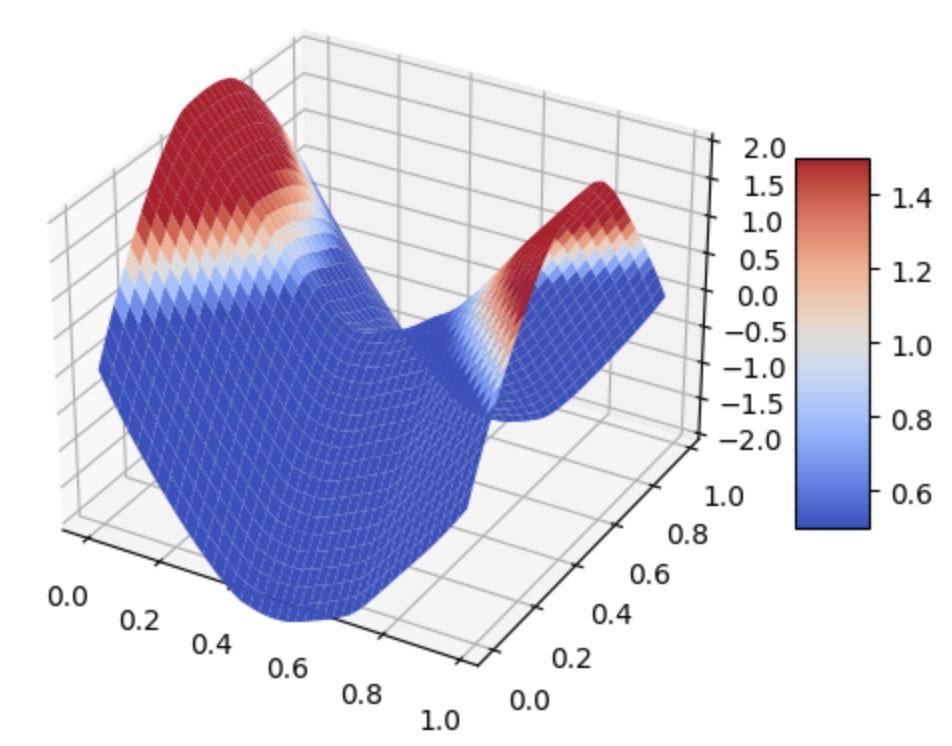}
    \caption{A Neural Network used to model a minimal surface for Scherk's first surface with 5000 training epochs.}
    \label{fig:scherknn5}
    \end{minipage}
    \hfill
    \begin{minipage}[b]{0.3\textwidth}
    \includegraphics[width = \textwidth, scale=0.2]{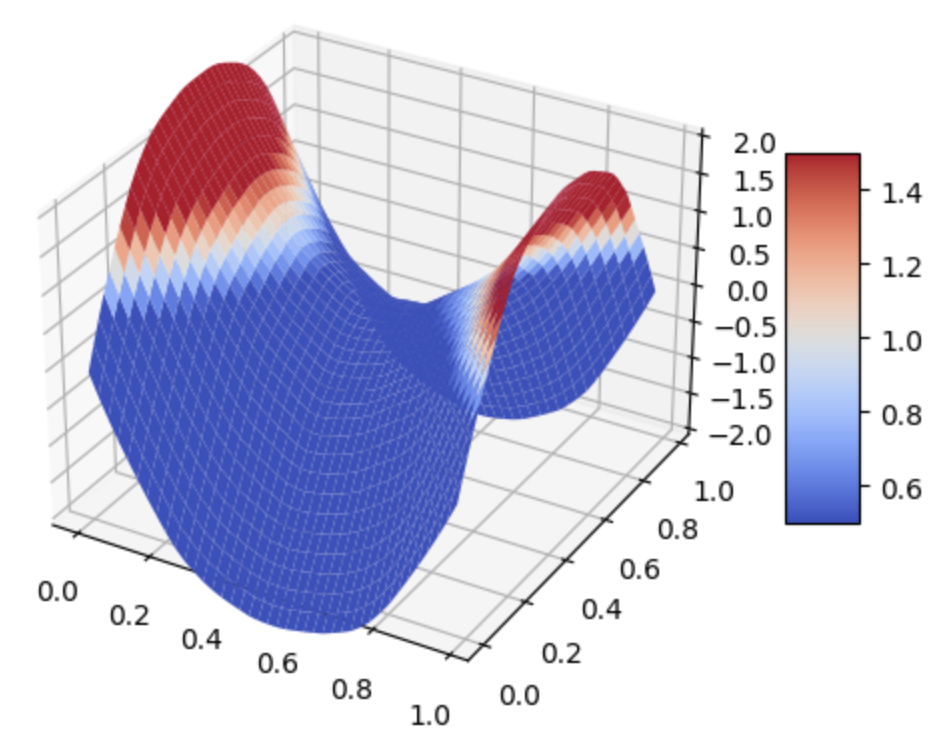}
    \caption{Another Neural Network for Scherk's first surface with 20,000 training epochs.} 
    \label{fig:scherknn20}
    \end{minipage}
    \hfill
    \begin{minipage}[b]{0.3\textwidth}
    \includegraphics[width = \textwidth, scale=0.2]{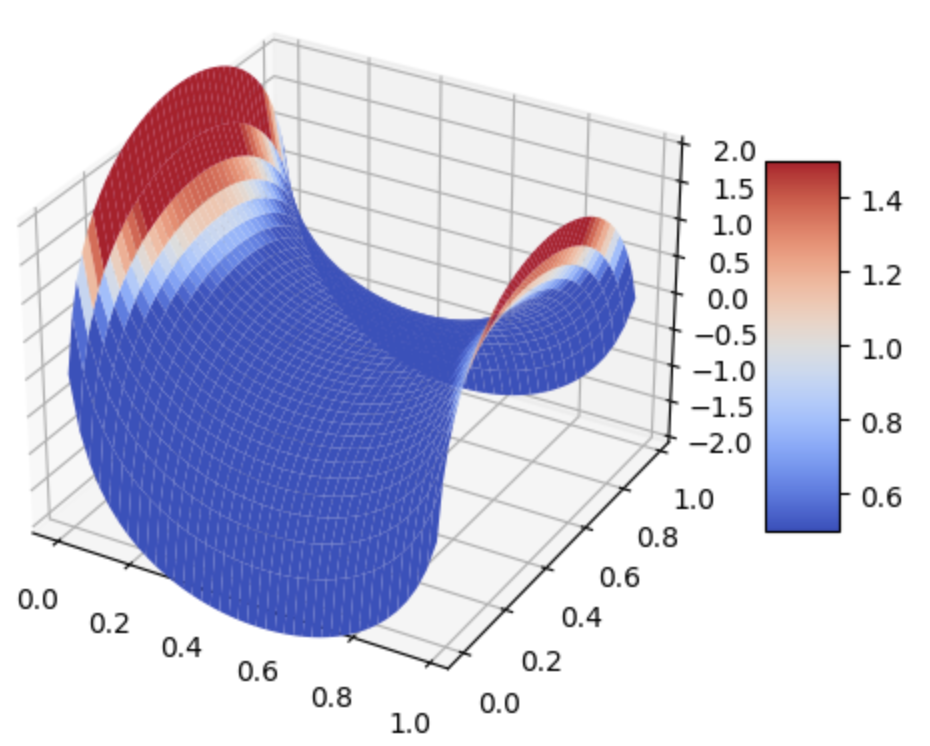}
    \caption{An original unit cell of Scherk's first minimal surface.} 
    \label{fig:scherkorg}
    
    \end{minipage}
        \centering
    \includegraphics[width = 0.8\textwidth]{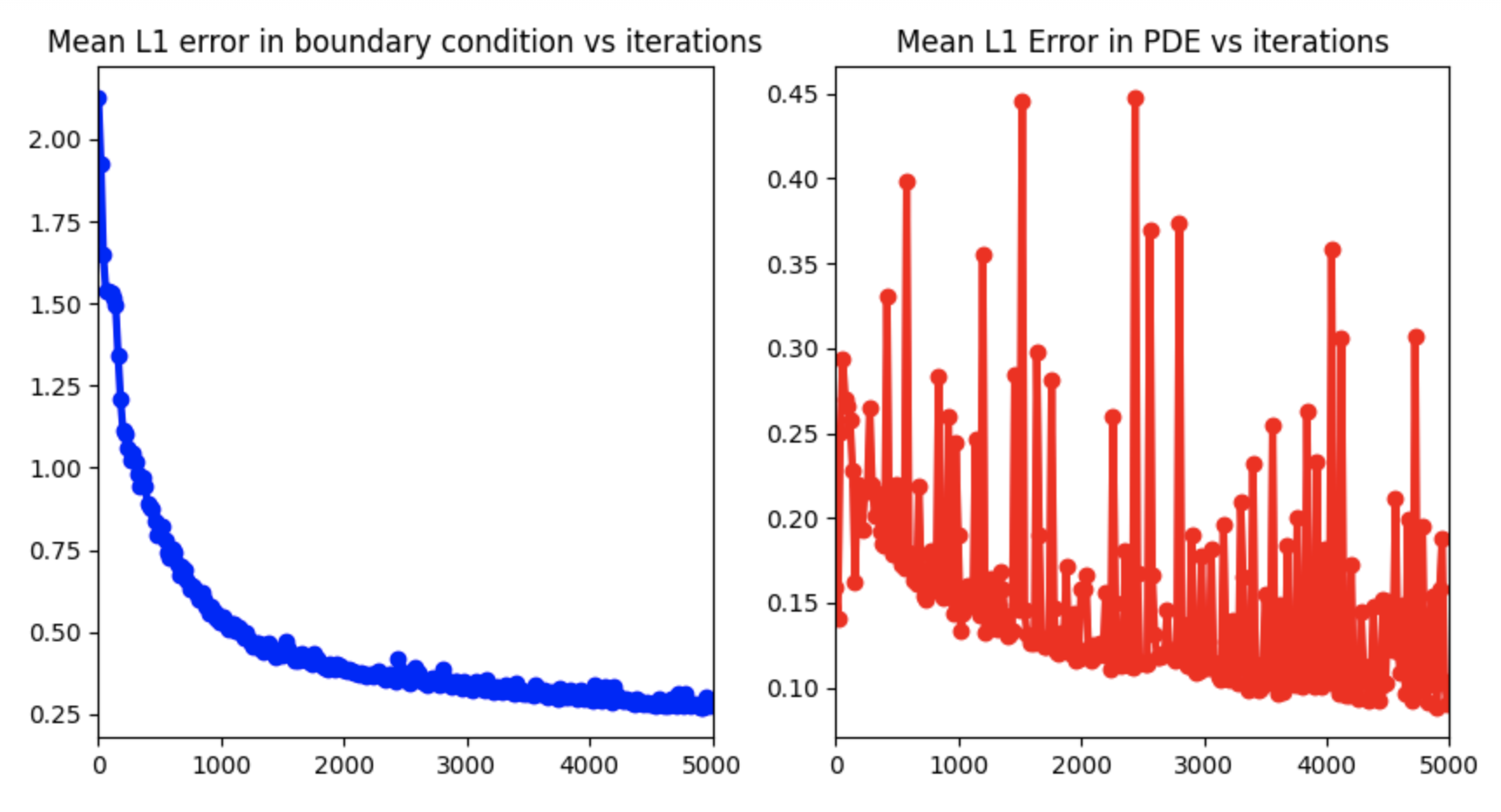}
    \caption{The unweighted boundary(left) and internal PDE(right) loss calculated for Scherk's minimal surface in figure \ref{fig:scherknn5} after 5000 epochs.}
    \label{fig:5loss}
    
    \centering
    \includegraphics[width = 0.8\textwidth]{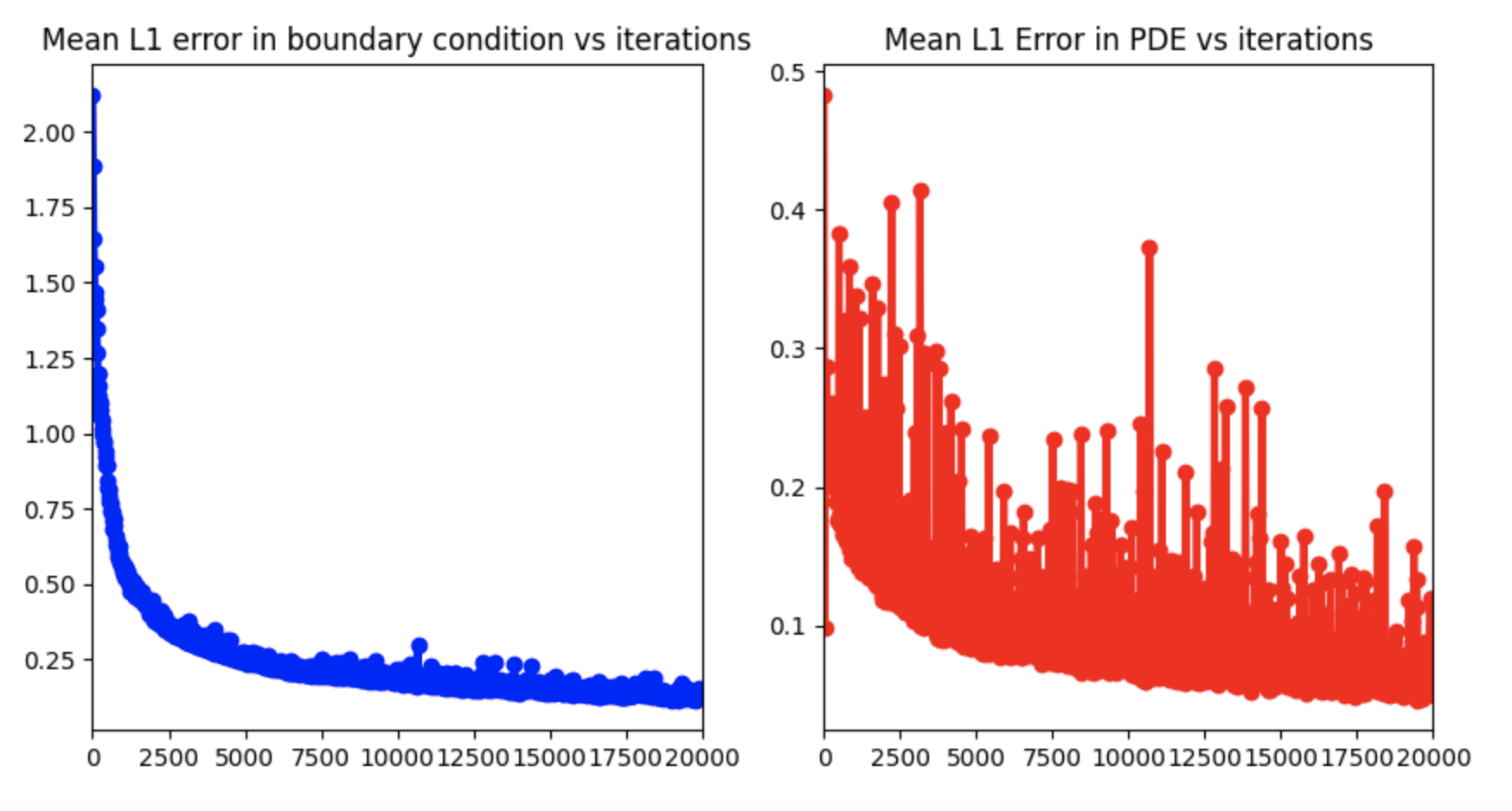}
    \caption{The unweighted boundary(left) and internal PDE(right) loss calculated for the same minimal surface but with 20,000 epochs.}
    \label{fig:20loss}
\end{figure}

In this paper we will introduce a total of seven examples of minimal surface plotting in 2D, 3D, and 4D domains. In all of these cases, we will be using the same optimizer and general training methods, with exceptions being a difference in the weight of the boundary loss versus the interior loss, optimizer learning rate, and number of epochs. 

To showcase the simplicity and potential complexity of our model, we will use the Pytorch library, a newer alternative to traditional libraries like Tensorflow. 

Our optimizer, as earlier mentioned, is Adam. As an upgraded version of the traditional Stochastic Gradient Descent (SGD), Adam delivers much more efficient results over a smaller amount of epochs by using adaptive moments. Another benefit is the constant micro-adjustments Adam makes to prevent overshooting the loss, avoiding SGD's potentially inaccurate results. 

All of the below examples were run on Deepnote's basic processing machine, with 5GB of RAM and a 2vCPU. No GPU was used, greatly slowing the training process. 

In each example, 1000 points have been picked from the interior of the ``frame'' to approximate $f(x)$ and 200 points are placed on each edge summing to a total of 800 boundary points in 2D, 2400 points in 3D, and 6400 points in 4D. 

It is important to note before moving on that while our input domains are in 2D, 3D, and 4D, all of the outputted minimal surfaces will be in 3D, 4D, and 5D respectively. 

\subsection{Two Dimensions}

\begin{figure}
    \centering
    \begin{minipage}[b]{0.45\textwidth}
    \includegraphics[scale=0.4, width = \textwidth]{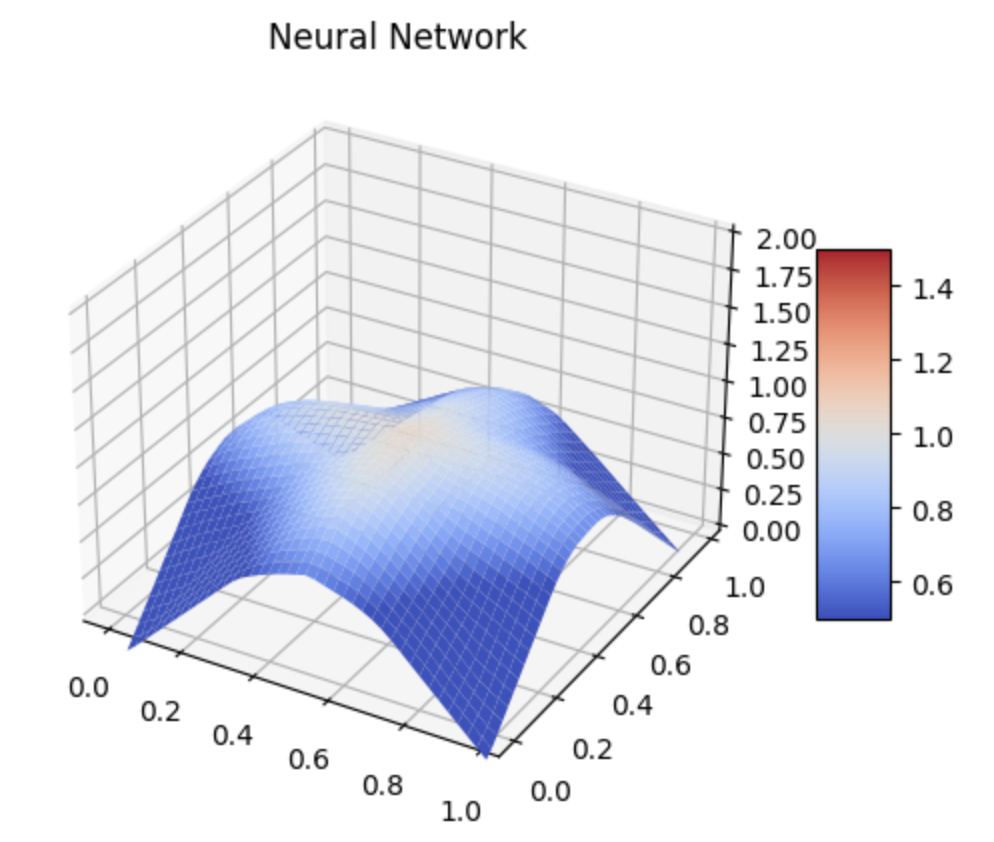}
    \caption{A Neural Network used to model a minimal surface with boundary condition \eqref{eq:firstbound} in a 2D input domain.}
    \label{fig:firstnn}
    \end{minipage}
    \hfill
    \begin{minipage}[b]{0.45\textwidth}
    \includegraphics[width = \textwidth, scale=0.2]{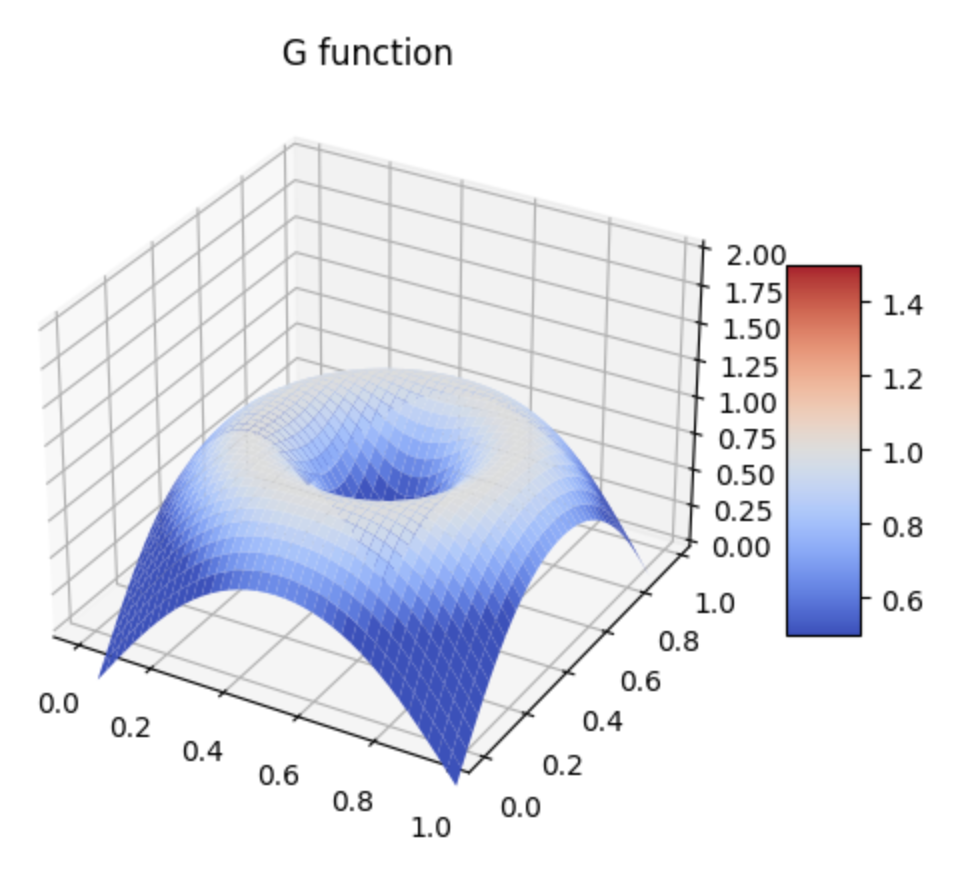}
    \caption{The original boundary condition \eqref{eq:firstbound} extended as a surface throughout the 2D domain}
    \label{fig:firstg}
    \end{minipage}
    \centering
    \includegraphics[width = \textwidth]{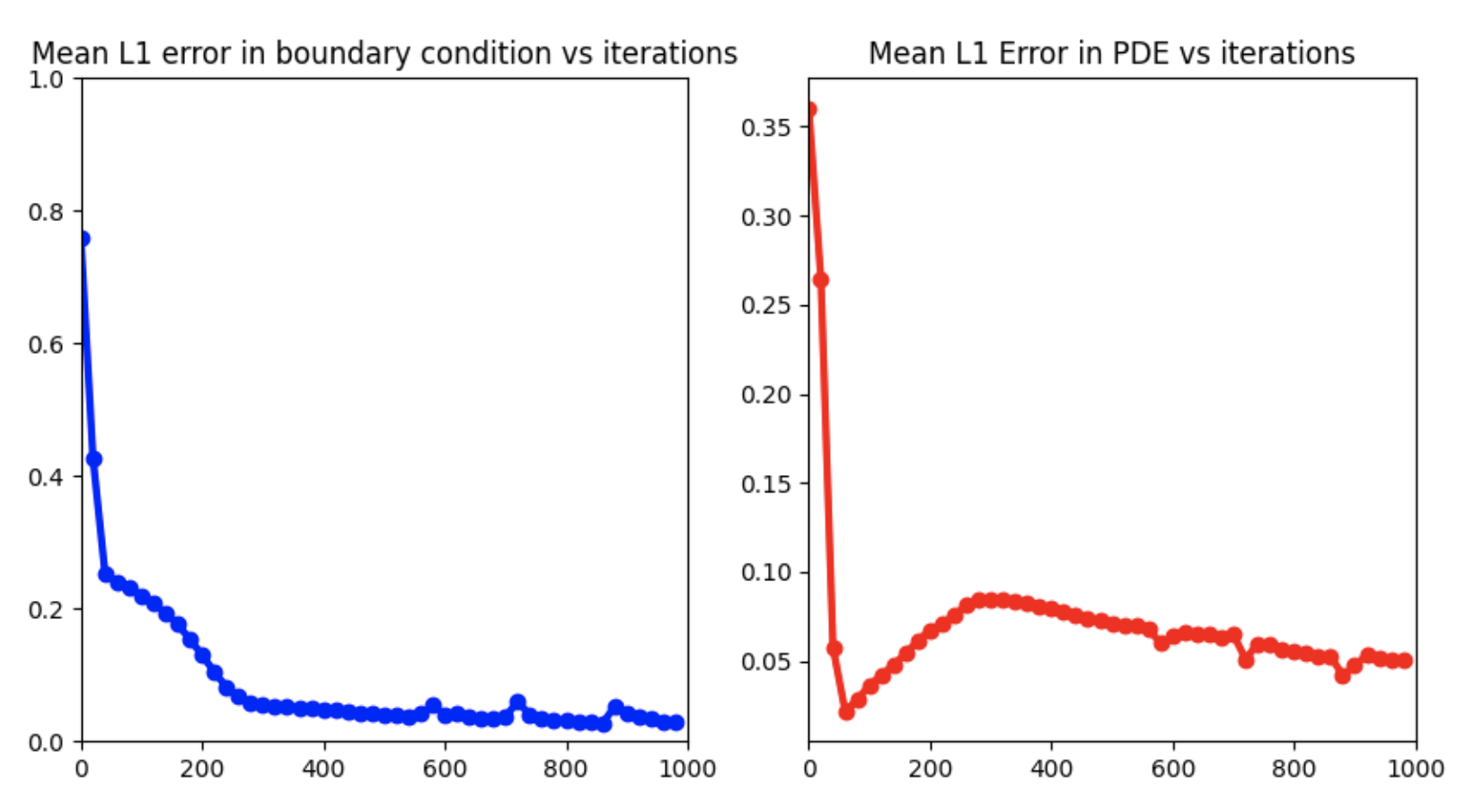}
    \caption{The unweighted boundary(left) and internal PDE(right) loss calculated for figure \ref{fig:firstnn} after 1000 epochs.}
    \label{fig:firstloss}
\end{figure}

\begin{figure}
    \centering
    \begin{minipage}[b]{0.45\textwidth}
    \includegraphics[scale=0.4, width = \textwidth]{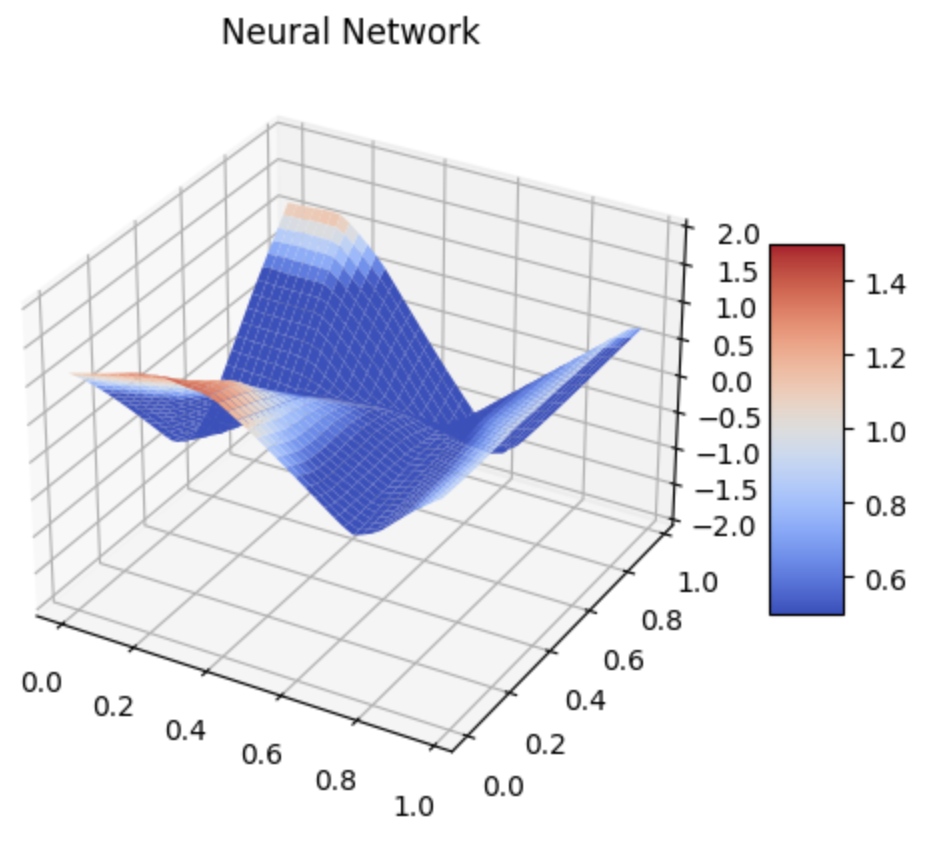}
    \caption{A Neural Network used to model a minimal surface with boundary condition equation \eqref{eq:secondbound} in a 2D input domain.}
    \label{fig:secondnn}
    \end{minipage}
    \hfill
    \begin{minipage}[b]{0.45\textwidth}
    \includegraphics[width = \textwidth, scale=0.2]{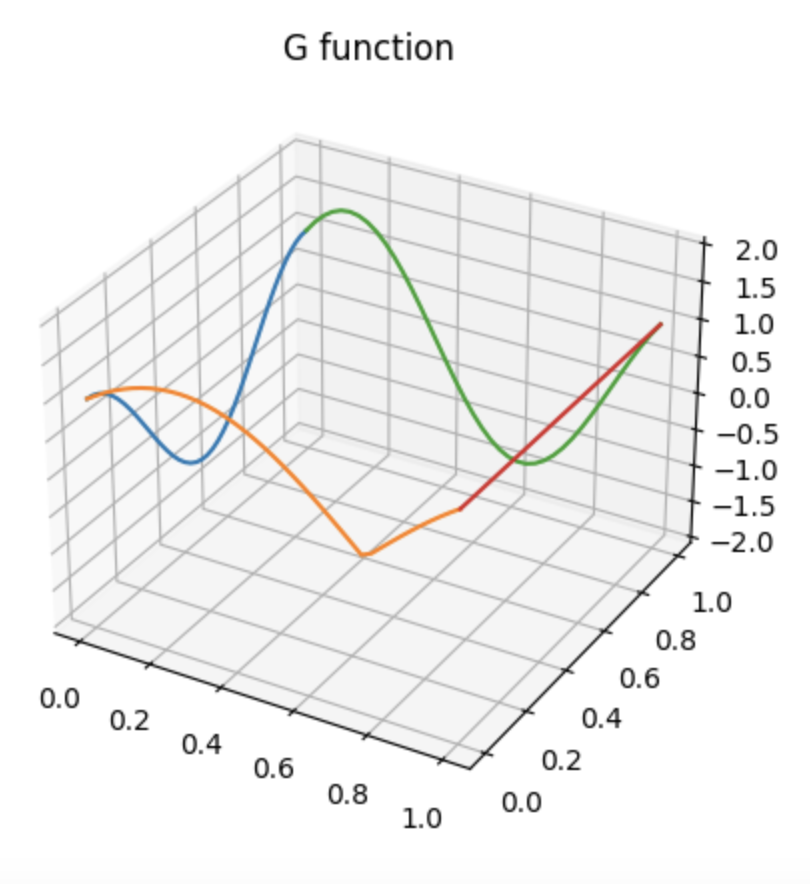}
    \caption{The original 4 boundaries graphed in the same 3D space. No coherent equation can be naturally defined to fill in the interior} 
    \label{fig:secondg}
    \end{minipage}
    \centering
    \includegraphics[width = \textwidth]{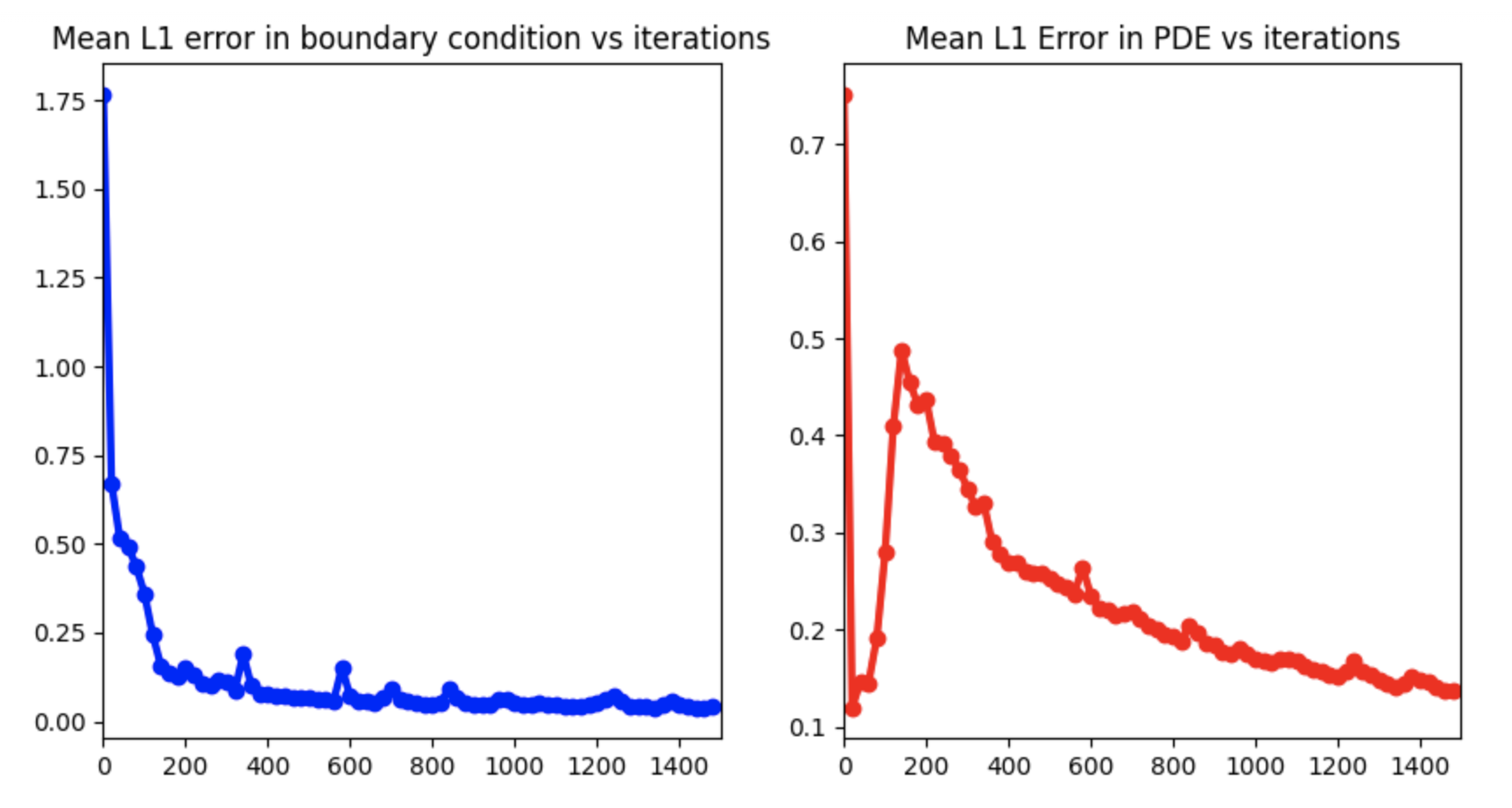}
    \caption{The unweighted boundary(left) and internal PDE(right) loss calculated for figure \ref{fig:secondnn} after 1500 epochs.}
    \label{fig:secondloss}
\end{figure}

Three examples of 2D will be used to showcase the basic usage of our programming and then apply it to a boundary condition where the minimal surface becomes apparent when compared to the boundary equations.

We will use Scherk's first minimal surface as an example to confirm that we are in fact, making an accurate minimal surface. The equation in question is

\begin{equation}
g(x) = \log{(\frac{\cos{x_{2}}}{\cos{x_1}})},  \quad f(x) = 0
\label{eq:scherk}
\end{equation}

This only shows one cell of the minimal surface, and in figures \ref{fig:scherknn5}, \ref{fig:scherknn20}, and \ref{fig:scherkorg}  the domain is selected as $[-1.5, 1.5]^2$ instead of $[-\frac{\pi}{2}, \frac{\pi}{2}]^2$ to prevent an indeterminate result. Both 5000 and 20,000 epochs of training are recorded for this equation at a learning rate of 0.001. For this scenario, the boundary weight is only $30\%$ of the PDE weight. 

The neural-network-generated minimal surface is visually identical to a single cell of Scherk's first surface.

Between the 5000 epoch, 20,000 epoch, and original Scherk surface, there is very little difference in the shapes of the graphs showing that our model has produced accurate results. Furthermore, the PDE loss in Figures \ref{fig:5loss} and \ref{fig:20loss} both converge, showing that the number of epochs is not a significant limiter on our accuracy excluding the spikes in loss caused by Adam. From this data, we can safely assume that using around 1000-2000 epochs for each future case will be sufficient.

Our second boundary condition is as follows
\begin{equation}
g(x) = \sin{ (\lVert5 (x-x_0)\rVert_{2})}, 
\quad x_{0} = (0.5, 0.5), \quad f(x) = 0
\label{eq:firstbound}
\end{equation} 
For this specific boundary example, our Adam optimizer has been set to a learning rate of 0.003. The MSE loss for the boundary condition will be weighed 3 times greater than the weight of the MSE loss for the internal minimal surface PDE. After 1000 iterations of neural network training, we are left with the minimal surface in Figure \ref{fig:firstnn}. 

Figure \ref{fig:firstg} is used as a foil to our surface-area conservative neural network. From this original equation we observe a significant difference between the expected value of the function should the entire surface follow the behavior of the boundaries alone and the minimal surface. 

Figure \ref{fig:firstloss} further consolidates our certainty of the accuracy of the model. As the number of epochs approaches the limit we have set for the current training session, the loss for the boundary condition, which we view as more important in not just the weight of the two MSE losses, but also integral to the fundamental definition of a "soap bubble" minimal surface, minimizes to a value below 0.04. With the average error calculated for every 20 epochs, we can see that in a short duration of training, about two minutes on a Macbook, even the fairly erratic internal PDE loss enters a decreasing trend that settles around 0.05.

Another example helps fully certify our neural network as an accurate method of analyzing different boundary conditions. 

The boundary conditions for this example are as follows, with y being a variable spanning from 0 to 1

\begin{equation}
\begin{aligned}
&x_{1} = 0, \quad g(x) = \sum_{i=1}^n \cos{(2\pi x_{i})}, \quad f(x) = 0\\
&x_{1} = 1, \quad g(x) = 1\\
&x_{2} = 0, \quad g(x) =\sum_{i=1}^n |(\cos{(\pi x_{i})} + \sin{(\pi x_{i})})|\\
&x_{2} = 1, \quad g(x) = \sin{(2\pi x_{1})} + \cos{(2\pi x_{1} })
\label{eq:secondbound}
\end{aligned}
\end{equation} 

Now we completely remove the center of our area and give each boundary a different equation to test the limits of the model. 

Our specific values of training plugged into the optimizer have been fine-tuned to provide a more accurate model of the minimal surface. Due to the nature of the boundary equations, we raise the weights to a 1:5 ratio in favor of the boundary.

Figure \ref{fig:secondnn} is the same as the earlier neural network; however, Figure \ref{fig:secondg} appropriately excludes any surface within the boundary frame. Referring back to the ``soap bubble'' analogy, if the boundaries of the 2D domain are the frames of a soap bubble mold, then the soap will follow the form presented in Figure \ref{fig:secondnn}. 

For this boundary condition, the model isn't able to quickly balance out the PDE loss within 1500 epochs resulting in an initial spike. But with the downward trend that is shown in the majority of the graph, we can be relatively certain that more epochs of training will balance out the loss. The need to optimize either the boundary or the internal PDE neglects the other parameter. The only method for creating a mutually inclusive solution that takes care of both necessities is to precisely weigh both of the MSE losses of each graph in a way that we as the researchers become satisfied with the overall minimization of both graphs. Through many cycles of trial and error, we have finally come to find some of the closest values to those that invoke the least overall loss.

With all of the added inconsistency in the frame, our model still holds up quite well in 2D situations, which ensures a great deal of practicality in potential real-world applications.

\subsection{Three Dimensions}

\begin{figure}
    \centering
    \begin{minipage}[b]{0.45\textwidth}
    \includegraphics[scale=0.4, width = \textwidth]{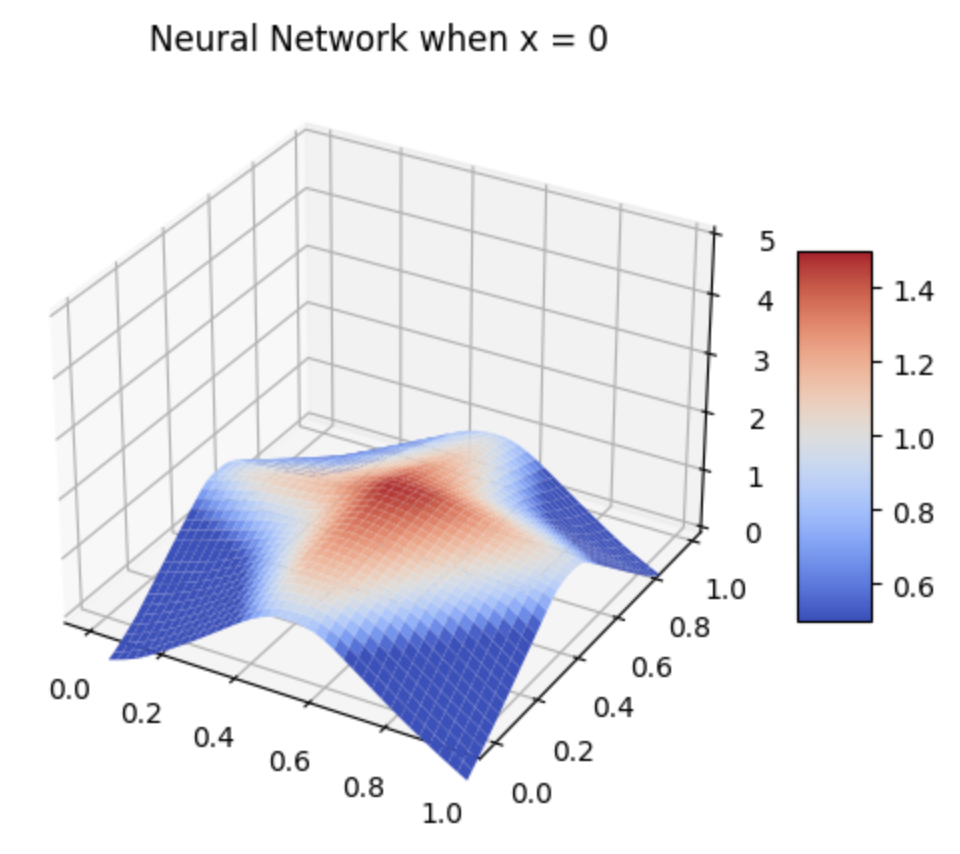}
    \caption{A Neural Network used to model a minimal surface with boundary condition \eqref{eq:thirdbound} in a 3D input domain. The face shown is when x is constantly 0.}
    \label{fig:thirdnn}
    \end{minipage}
    \hfill
    \begin{minipage}[b]{0.45\textwidth}
    \includegraphics[width = \textwidth, scale=0.2]{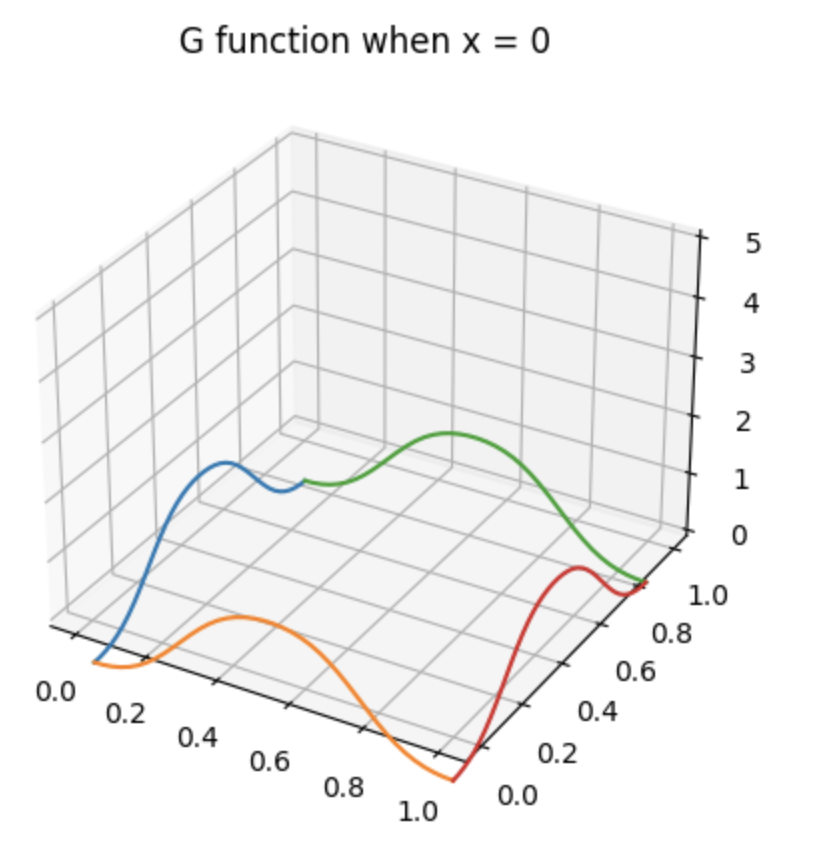}
    \caption{The original boundary condition extended as a surface throughout the 3D domain. The face shown is corresponding to the one shown in the neural network graph.}
    \label{fig:thirdg}
    \end{minipage}
    \centering
    \includegraphics[width = \textwidth]{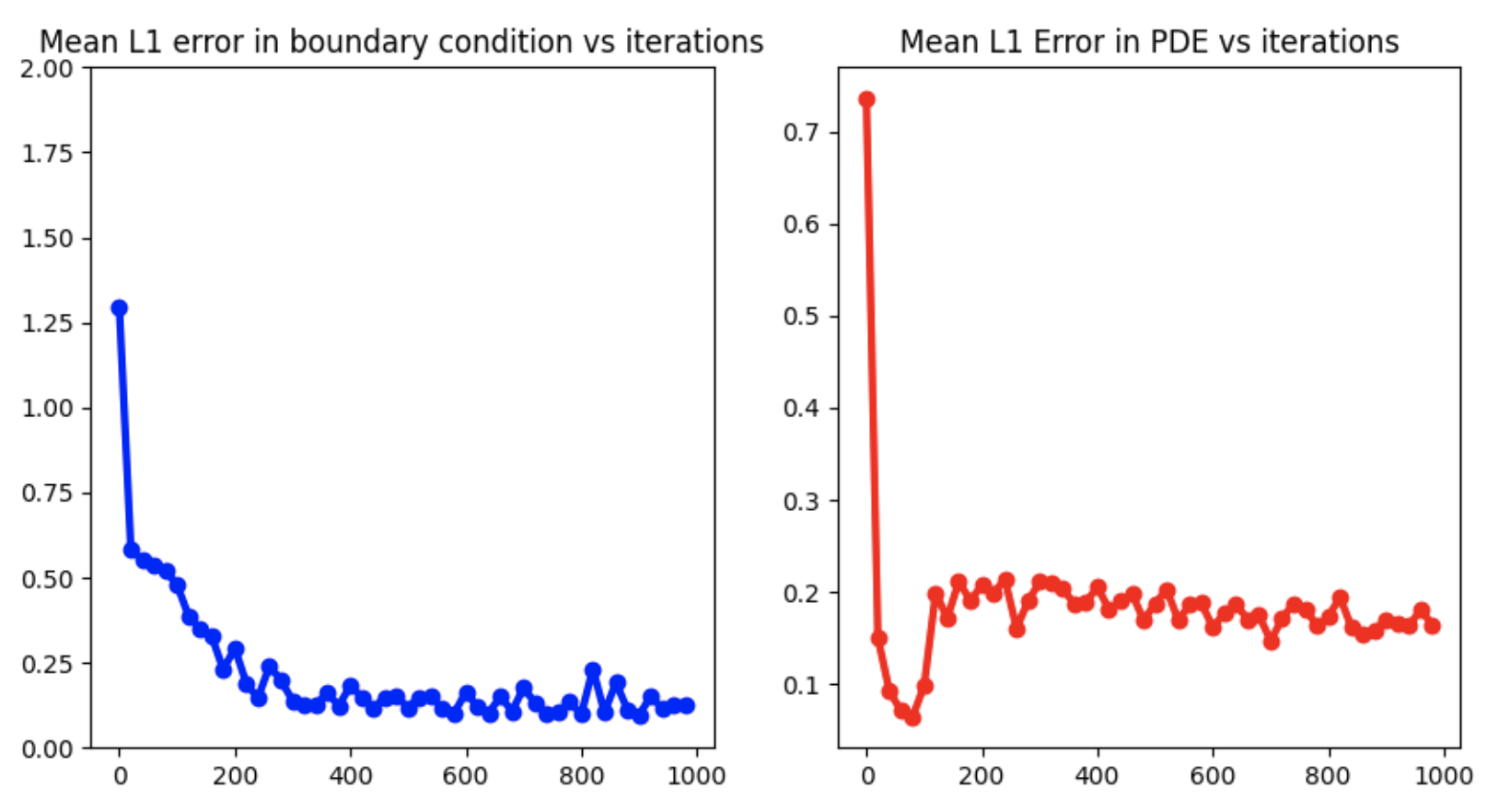}
    \caption{The unweighted boundary(left) and internal PDE(right) loss calculated for \ref{fig:thirdnn} after 1000 epochs.}
    \label{fig:thirdloss}
\end{figure}

Once our boundary condition reaches 3 dimensions, the method of visualizing the domain is no longer as intuitive. Instead of having a four-sided square as a base domain with boundaries that go from 0 to 1 in $x$ and $y$ axis and a function output reaching into the $z$-axis, we now have to build a cube-like structure for each boundary equation. Each edge of a cube with 3 coordinates each being either 0 or 1 is separately applied to our boundary condition defined in equation \eqref{eq:thirdbound}. Every boundary equation has been carefully set so that the boundary values will be equal at each corner to prevent decomposition of the cubical shape which can be seen as a leak where the frame becomes disconnected.
\begin{equation}
g(x) = \sin{(\cos{(2\pi \sum_{i=1}^{3} |x_i-x_{0i}|)})}, 
\quad x_{0} = (0.5, 0.5, 0.5), \quad f(x) = 0.
\label{eq:thirdbound}
\end{equation} 

Similarly to the two-dimensional examples, we will use largely the same machine learning settings. This time, our optimizer will be Adam with a learning rate of 0.001 over the course of 1000 epochs. We have also given the two values the exact same weight. 

\begin{figure}
    \centering
    \begin{minipage}[b]{0.45\textwidth}
    \includegraphics[scale=0.4, width = \textwidth]{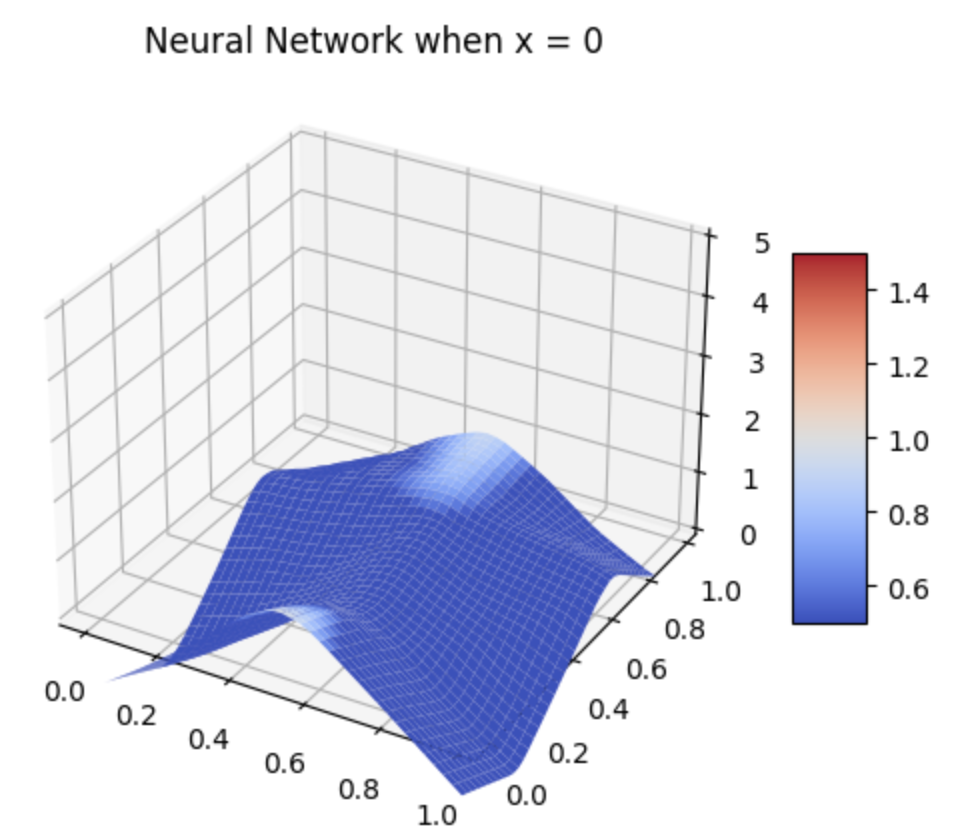}
    \caption{A Neural Network used to model a minimal surface with boundary equation \eqref{eq:fifthbound} in a 3D input domain. The current face shown is that when variable x is set to 0 while y and z are variables}
    \label{fig:fifthnnv1}
    \end{minipage}
    \hfill
    \begin{minipage}[b]{0.45\textwidth}
    \includegraphics[width = \textwidth, scale=0.2]{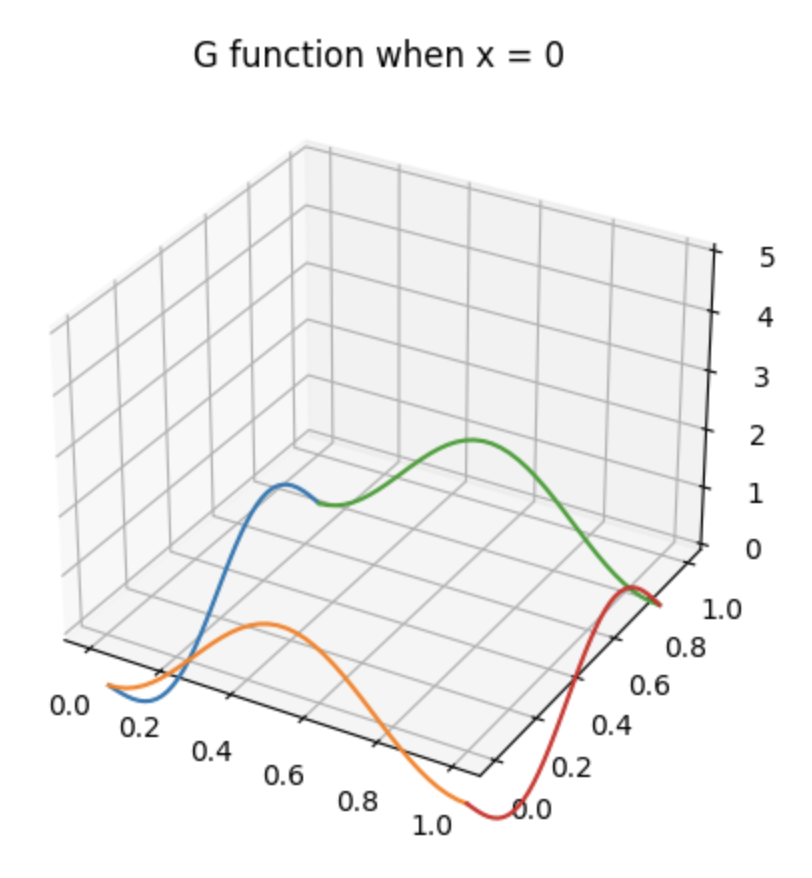}
    \caption{The original boundary condition extended as a surface throughout the 3D domain. The face shown is when x is the constant 0, the same as in the neural network to the left.}
    \label{fig:fifthgv1}
    \end{minipage}
    \centering

    \centering
    \begin{minipage}[b]{0.45\textwidth}
    \includegraphics[scale=0.4, width = \textwidth]{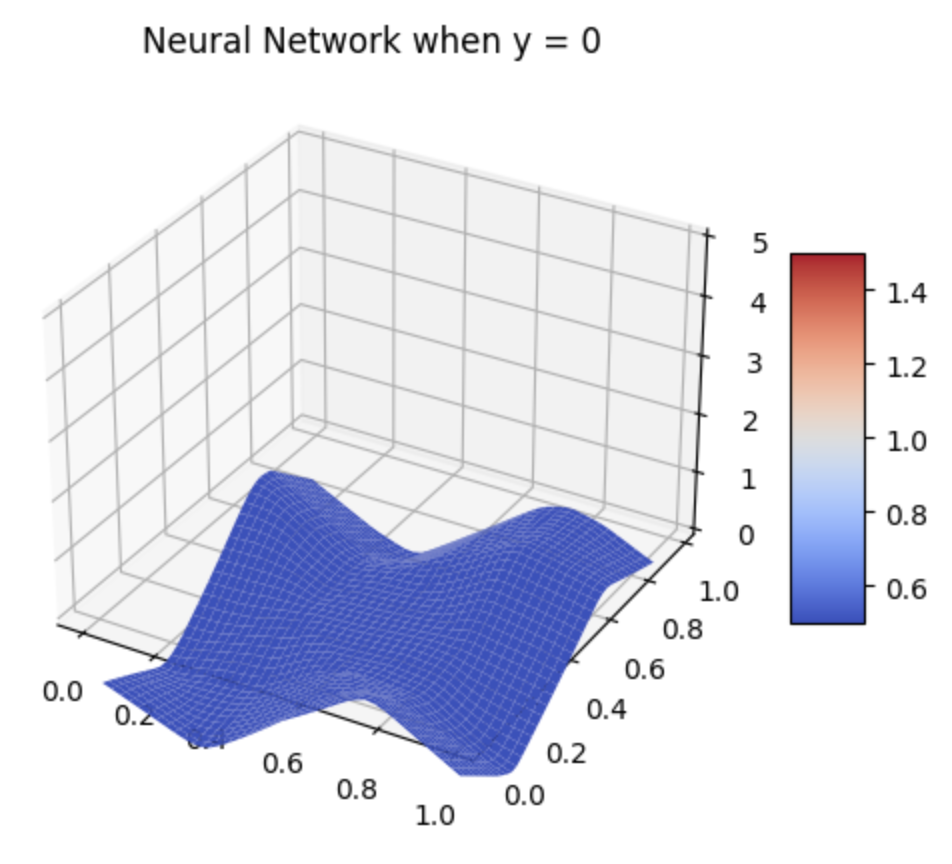}
    \caption{A Neural Network used to model a minimal surface with boundary equation \eqref{eq:fifthbound} in a 3D input domain. Here, y is held constant as 0 while x and z are variable.}
    \label{fig:fifthnnv2}
    \end{minipage}
    \hfill
    \begin{minipage}[b]{0.45\textwidth}
    \includegraphics[width = \textwidth, scale=0.2]{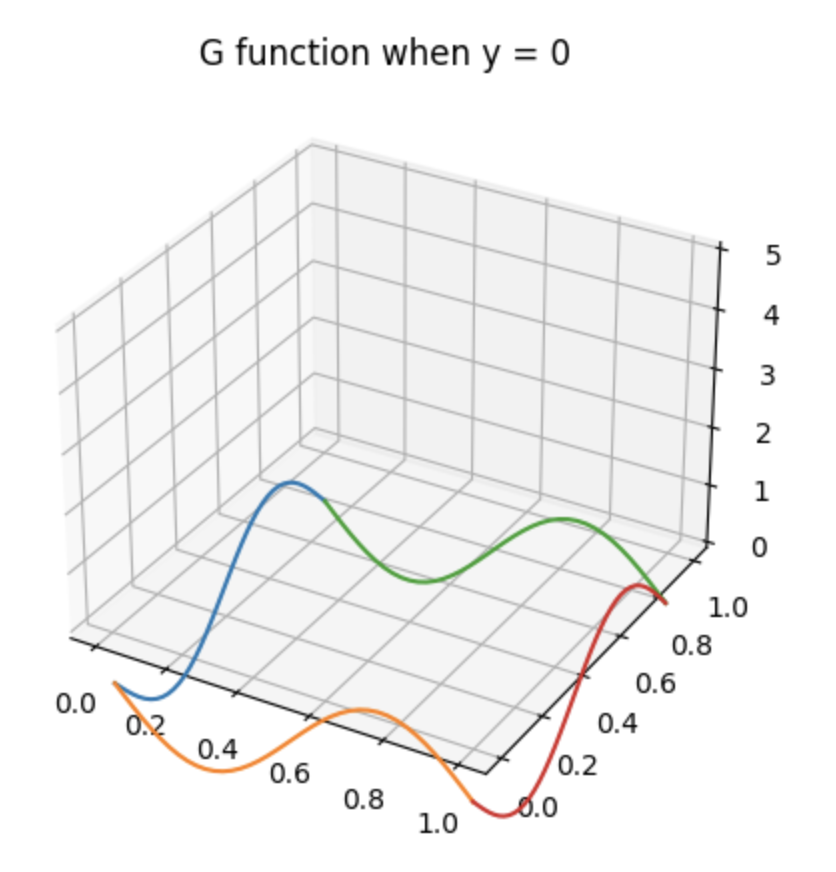}
    \caption{The original boundary condition extended as a surface throughout the 3D domain. The face shown is the same as the neural network on the left, when y = 0.}
    \label{fig:fifthgv2}
    \end{minipage}
    \centering
\end{figure}

\begin{figure}
    \centering
    \begin{minipage}[b]{0.45\textwidth}
    \includegraphics[scale=0.4, width = \textwidth]{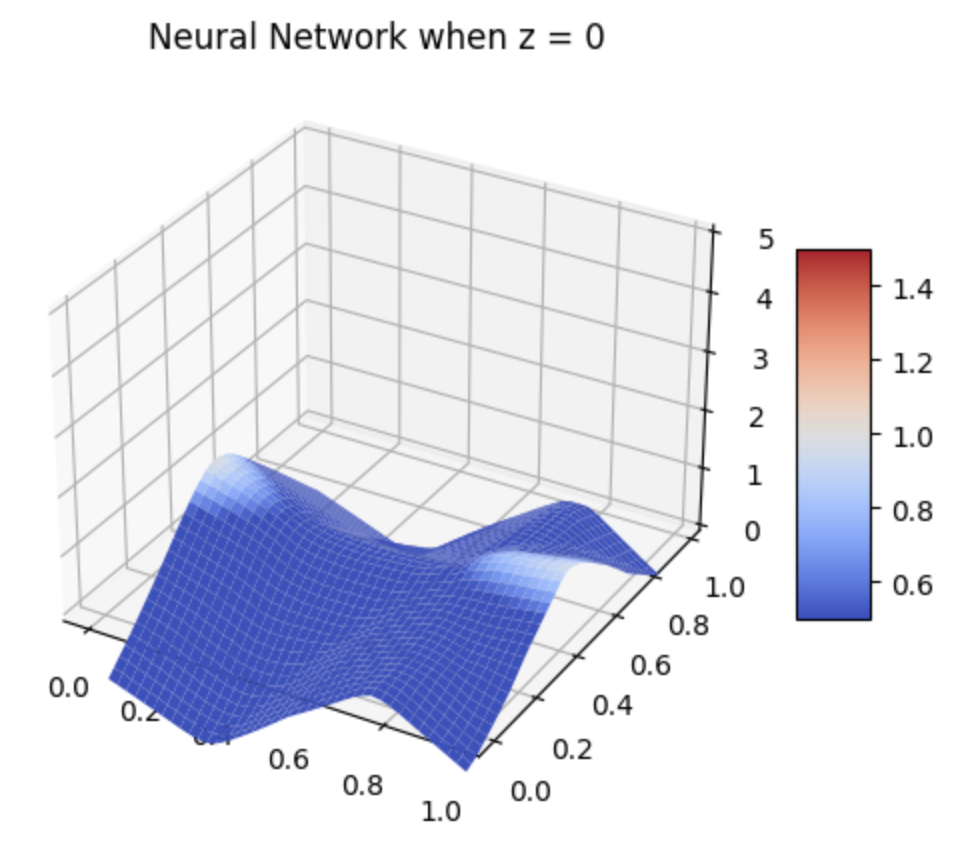}
    \caption{A Neural Network used to model a minimal surface with boundary equation \eqref{eq:fifthbound} in a 3D input domain. Z is held constant with value 0 while x and y are variable.}
    \label{fig:fifthnnv3}
    \end{minipage}
    \hfill
    \begin{minipage}[b]{0.45\textwidth}
    \includegraphics[width = \textwidth, scale=0.2]{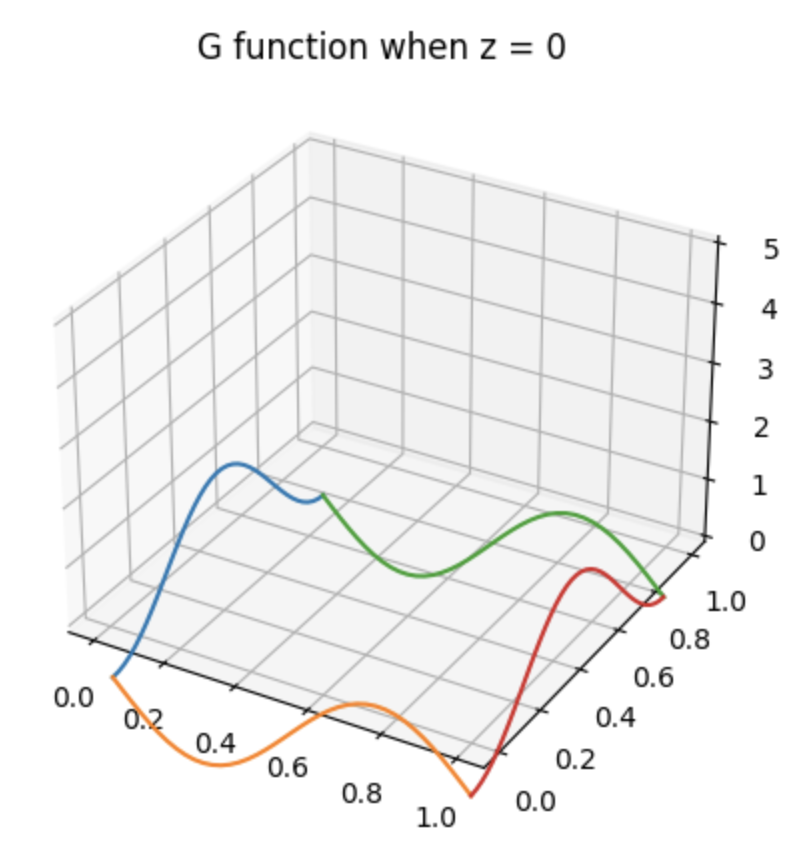}
    \caption{The original boundary condition extended as a surface throughout the 3D domain and projected onto a 2D graph domain for a constant z variable set to 0.}
    \label{fig:fifthgv3}
    \end{minipage}
    \centering

    \includegraphics[width = \textwidth]{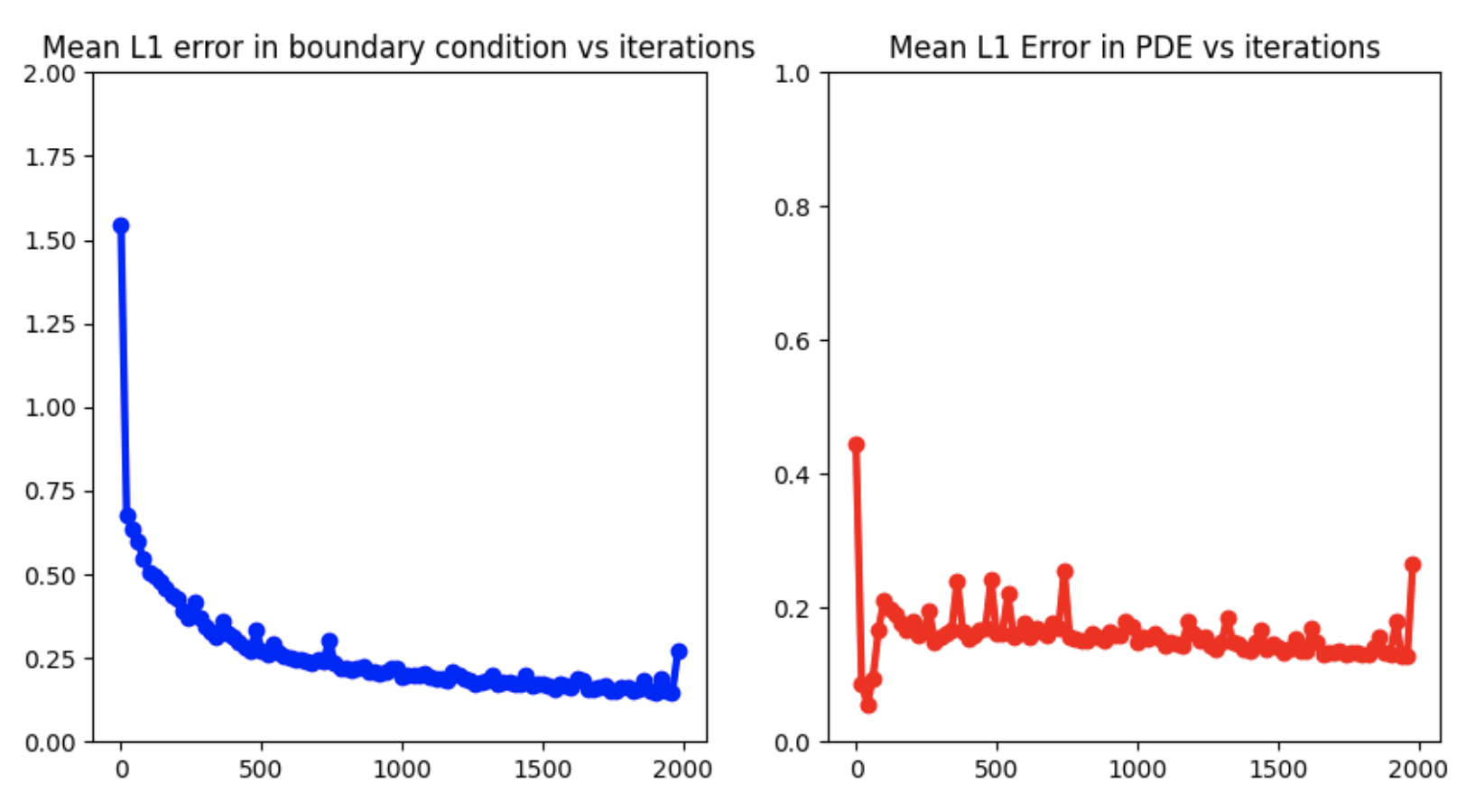}
    \caption{The unweighted boundary(left) and internal PDE(right) loss calculated for the neural networks in figures \ref{fig:fifthnnv1}, \ref{fig:fifthnnv2}, and \ref{fig:fifthnnv3} after 2000 epochs.}
    \label{fig:viewloss}
\end{figure}
\begin{figure}
    \centering
    \begin{minipage}[b]{0.45\textwidth}
    \includegraphics[scale=0.4, width = \textwidth]{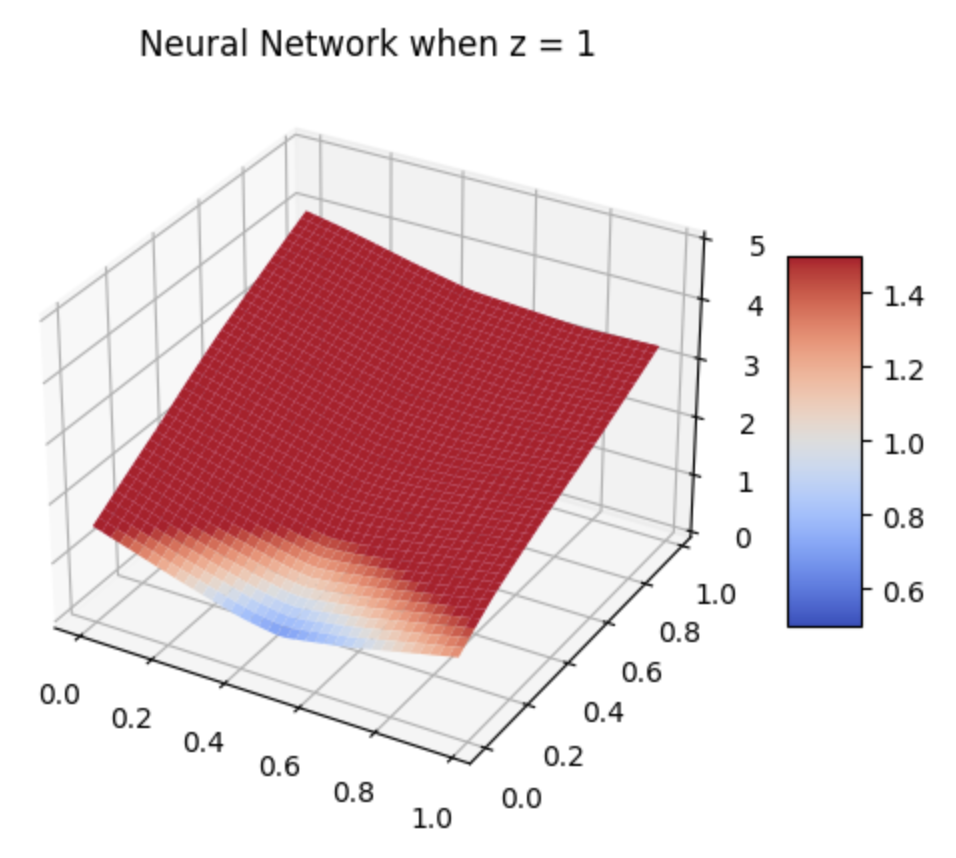}
    \caption{One face of a Neural Network where z = 1. The bottom-left edge is lacking in details and does not resemble the intended boundary.}
    \label{fig:seventhnn}
    \end{minipage}
    \hfill
    \begin{minipage}[b]{0.45\textwidth}
    \includegraphics[width = \textwidth, scale=0.2]{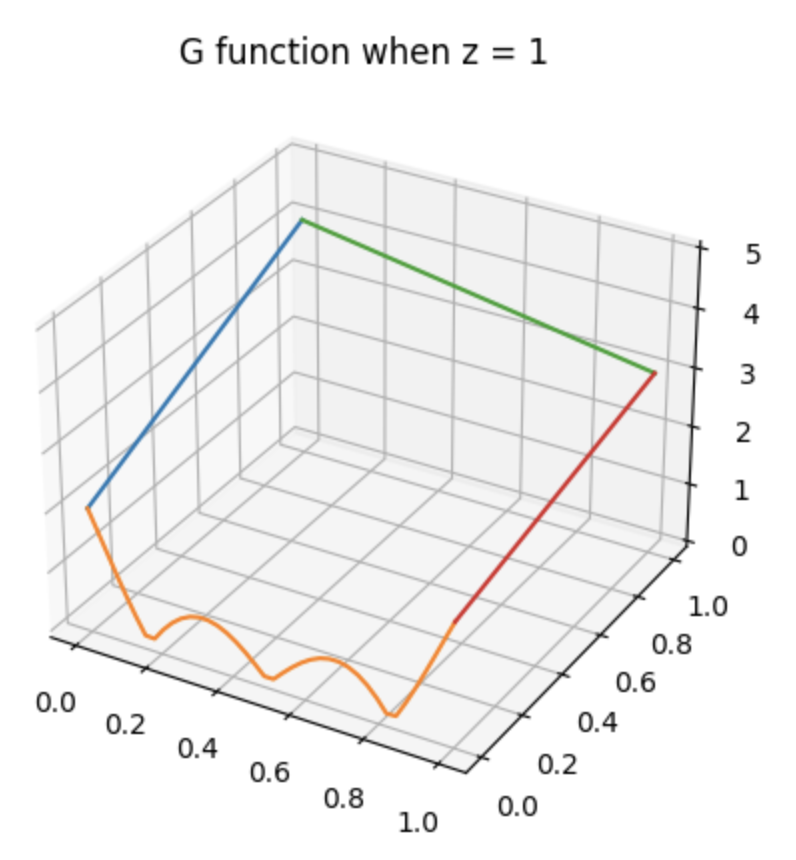}
    \caption{The original boundary condition as 4 different curves in a 3D space, corresponding to the face of the minimal surface where z = 1.}
    \label{fig:seventhg}
    \end{minipage}
\end{figure}

One of the greatest limitations of this strategy is that the output $g(x)$ enters the fourth dimension, one that we cannot visualize as humans and therefore cannot showcase in a single graph. The simplest method would be to adopt a strategy similar to taking partials of any multi-variable equation. We ``fix'' certain variables to constants and in that position, we find the 3D graph with a 2D variable input domain. 

In Figure \ref{fig:thirdg} we plot our initial g function that has a 3-variable input. We can fix x to be zero and produce a similar minimal surface graph as in the previous section. Note that the partial surfaces in Figures \ref{fig:thirdnn} and \ref{fig:thirdg} present only one of the 6 faces of the "cube" of boundary conditions and the z-axis in both of them is actually the fourth dimension that we cannot visualize. 

An extended application of our model defines a boundary equation that has differently shaped minimal surfaces on our 6 faces of the final minimal surface. Each variable of the three is synthesized into the equation 

\begin{equation}
g(x) = \sin{(2\pi x_1)} + \cos{(2\pi x_2)} + \sin{(2\pi x_3)}
, \quad f(x) = 0
\label{eq:fifthbound}
\end{equation} 

 When $x_1$, $x_2$, or $x_3$ is set to 0 or 1, the resulting aspect of sine or cosine is instantaneously equal to 0, with the combination of the different trigonometric functions for different variables creating the scenes in the different faces of Figures \ref{fig:fifthgv1}, \ref{fig:fifthgv2}, and \ref{fig:fifthgv3}. 

 As the increased demands for our PINN strain the ability to approximate the different surfaces, we change the dimensions of the optimizer accordingly. Our $g(x)$ loss weight is changed to a new value that carefully balances both errors without completely fulfilling either: a 1:1.5 ratio between PDE weight and $g(x)$ weight. Furthermore, the original 1000 iterations are now doubled to 2000, increasing the time required to fully train the model to over 6 minutes. 

 Notice that the differences in Figures \ref{fig:fifthnnv1}, \ref{fig:fifthnnv2}, \ref{fig:fifthnnv3} and their corresponding boundary equations are quite obvious when the original boundary function begins to pursue overly-exaggerated protrusions. While our minimal surface looks quite ``ugly'' compared to the boundary surfaces, there are small sections in the internal figure that indicate where area has been conserved.

 With these two examples we can build a basic understanding of each face of the higher-dimension minimal surface, but as our knowledge remains limited, all approximations or predictions we make are only speculative. 
 
 There are still some other limitations to using PINN, as are shown in Figures \ref{fig:seventhnn} and \ref{fig:seventhg}. Here, one edge of the 3D surface cannot be maintained after many epochs of training. This phenomenon can be observed due to the number of edges that have to be accounted for. Now that the model has a 3D bubble frame, the 12 edges comprising said frame each vie for importance during backpropagation. In that case, unnatural details will be mostly ignored by the PINN. 

 Despite many tests to overcome this limitation, such as raising the boundary weight to 10, 100, and even 10,000 times greater than the internal weight as well as decreasing the learning rate to 0.0001, the model still is unable to resolve the differences between the given G-function and neural network. 

 Another limitation not shown involves disconnects between the boundaries of different frames. Should the values on any corner of the frames not match up between the multiple equations at that position, the PINN will likewise struggle to produce accurate results.
 
The entire figure cannot be as simply defined as attaching each of our minimal surfaces to the faces of the cube and extruding into the 3 dimensions we currently have, but will require more advanced methods of visualizing the higher dimensions to truly understand. Nevertheless, the model's applications can still continue in this fashion.

\subsection{Four Dimensions}
 \begin{figure}
    \centering
    \begin{minipage}[b]{0.45\textwidth}
    \includegraphics[scale=0.4, width = \textwidth]{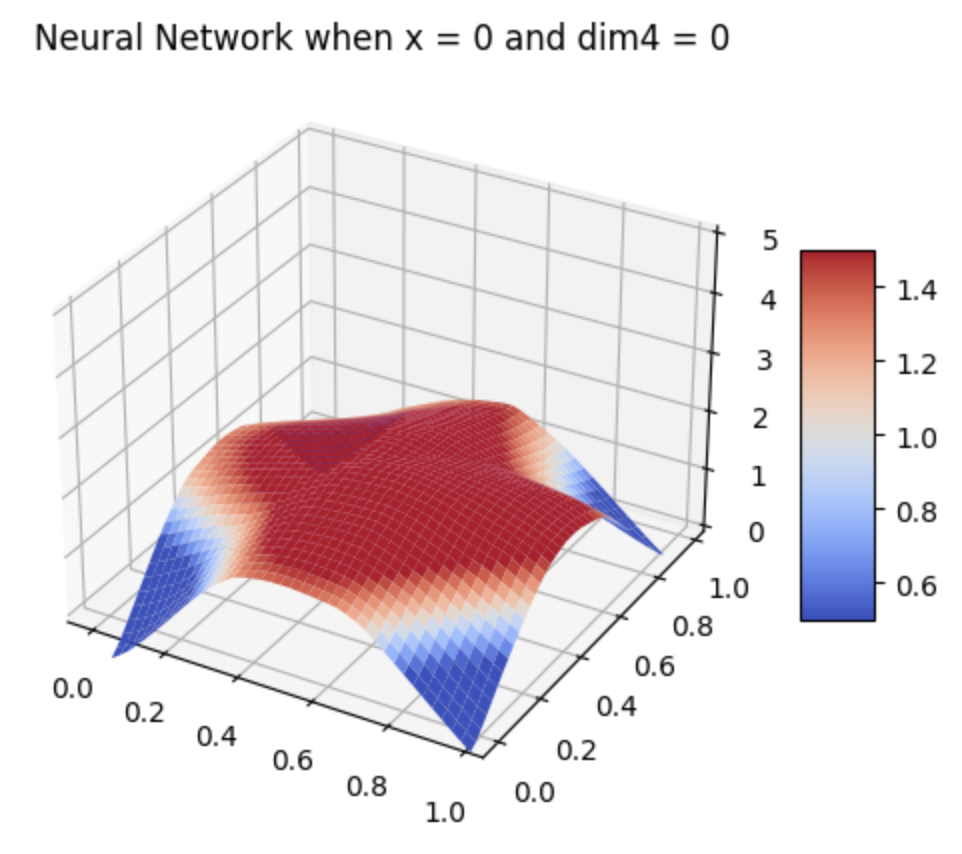}
    \caption{A Neural Network used to model a minimal surface with boundary condition \eqref{eq:fourthbound} in a 4D input domain. In this example, x and the 4th dimension are both set to 0.}
    \label{fig:fourthnn}
    \end{minipage}
    \hfill
    \begin{minipage}[b]{0.45\textwidth}
    \includegraphics[width = \textwidth, scale=0.2]{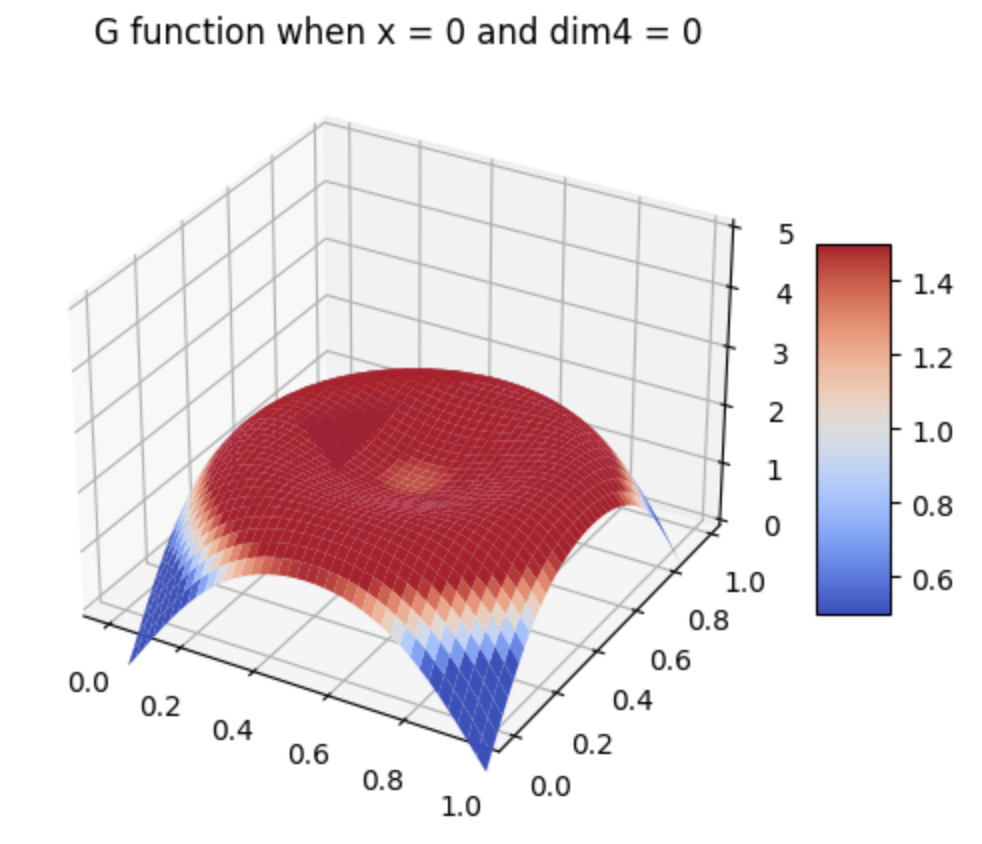}
    \caption{The original boundary condition extended as a surface throughout the 4D domain}
    \label{fig:fourthg}
    \end{minipage}
    \centering
    \includegraphics[width = \textwidth]{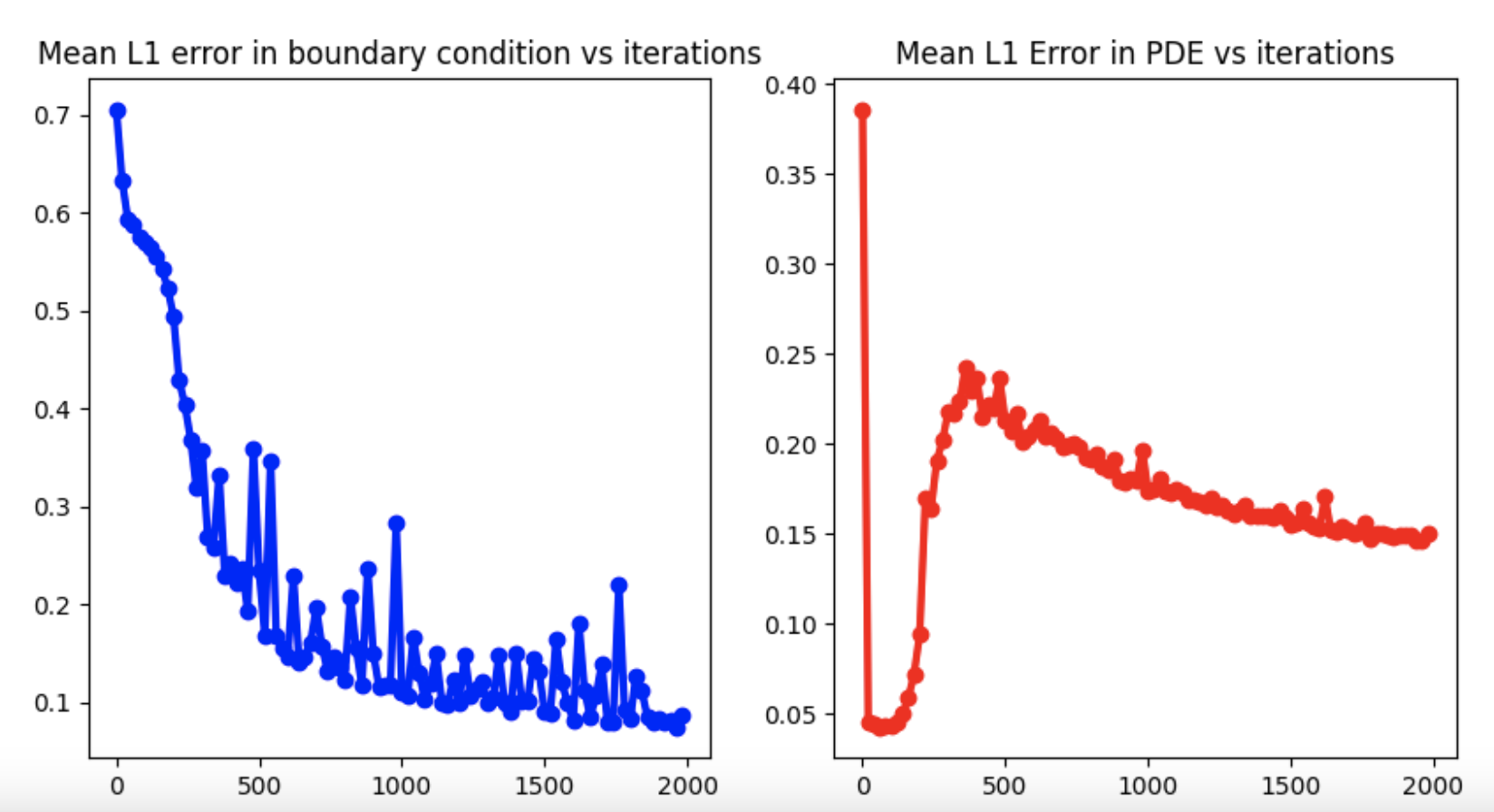}
    \caption{The unweighted boundary(left) and internal PDE(right) loss calculated for \ref{fig:fourthnn} after 1000 epochs.}
    \label{fig:fourthloss}
\end{figure}

To finish off our exploration into minimal surfaces, an input in 4 dimensions is the largest input that we attempt to model in this paper. Anything above that would begin to exceed the computing power of a Macbook. While we would be successful in such endeavors, it remains more efficient to showcase our model in such a way as to leave room for further improvement instead of blindly complicating without a true basis in the understanding of what our graphs would show. It is already hard enough to visualize the 4th dimension domain reaching into the 5th dimension, but we can still visualize some of the soap bubble. The strategy of fixing variables will continue to be our method of showcasing our results. 

Just as before, we will be building a frame in 4D for our boundary condition. Since we cannot visualize what the boundary will look like, we can only imagine it as extruding a cube into the 4th dimension, just as when the 2D square frame turned into the 3D cube frame.

Our last equation will take the form of 
\begin{equation}
g(x) = 2\sin{ (10 \lVert(x-x_0)\rVert_{2})}, 
\quad x_{0} = (0.5, 0.5, 0.5, 0.5), \quad f(x) = 0
\label{eq:fourthbound}
\end{equation} 
as an ode to our original 2-dimensional boundary equation. Having 2 more variables and changing one of the constants has changed the appearance and range of the output surface, but the continuity between the two models can still be seen. The boundary is weighed 3 times greater than the internal PDE and the Adam optimizer is set to 0.001 learning rate for 2000 epochs.

It can also be noted that the indent in Figure \ref{fig:fourthg} is completely gone in Figure \ref{fig:fourthnn}. The entire figure has also been flattened out in order to save surface volume around the highly-curving corner areas. The subtle differences in between the right figure, which is already quite smooth, and the left indicate the preciseness of which our model can seek out small imperfections and train them until they form a semblance of what we want. 

\begin{table}[]
    \centering
    \begin{tabular}{c|c|c|c}
         Example & Learning Rate & Weight Ratio & Epochs\\
         \hline
         2D - 1 & 0.001 & 0.3/1 & 500/20000\\
         2D - 2 & 0.003 & 1/3 & 1000\\
         2D - 3 & 0.003 & 1/5 & 1500\\
         3D - 1 & 0.001 & 1/5 & 1000\\
         3D - 2 & 0.001 & 1/1.5 & 2000\\
         3D - 3 & 0.0001 & 1/1 & 1000\\
         4D - 1 & 0.001 & 1/3 & 2000
    \end{tabular}
    \caption{List of Hyperparameters}
    \label{tab:para}
\end{table}

\section{Summary}

PINN is one of the most effective methods of solving the minimal surface equation, as has been seen through our multiple examples. Its computational efficiency and accurate results make it a strong contender for future PDE research. 

By carefully controlling the balance between the importance of weights, we have created a model that can intake results from a given domain frame and extend it through a small area in 2-4 dimensions. The boundaries of said frame can be given a variety of different equations governing how they enter the next dimension, and the model will in turn quickly and accurately generate the minimal surface.

In this work, we have described multiple examples of calculating minimal surface equations in spite of non-optimized boundary frames projecting each example as a convoluted shape into the next dimension. Our examples compound on the difficulty of their predecessors to provide deeper glances into the higher dimensions. By using minimal surfaces, we can see what soap bubbles in higher dimensions would look like in our 3 dimensions through snippets drawn where certain variables are fixed to constants. 

The overall performance of our methods still lacks robustness in high-stakes situations, so we will need newer methods for better results. Whether it is increasing the size of the neural network involved or redesigning the boundary loss equation, there are many ways to improve upon our base understanding of higher-dimensional minimal surfaces. Such methods still need to be developed, and with those advancements, we may gain better insights into the 4th dimension and beyond.

We hope that our work can be continued on with more optimized methods of visualizing the PINN such that audiences will finally visualize the higher dimension outputs without being limited to only seeing part of the figure.

\bibliographystyle{abbrv}
\bibliography{proj_ref}
\end{document}